\documentclass[11pt,a4paper]{amsart}
\pdfoutput=1
\usepackage[utf8]{inputenc}
\usepackage[T1]{fontenc}
\usepackage{amssymb,amsthm,amsmath}
\usepackage{tikz}

\usetikzlibrary{shapes.geometric, arrows}
\usepackage{url}
\usepackage{hyperref}
\usepackage[a4paper]{geometry}
\newtheorem{theorem}{Theorem}[section]

\newtheorem{contraction}[theorem]{Contraction Rules}
\theoremstyle{contraction}

\DeclareMathOperator{\sgn}{sgn}
\usepackage{pdfsync}

\begin{document}

\title{Mating of discrete trees and walks in the quarter-plane}
\author{Philippe  Biane}
\email  {biane@univ-mlv.fr }
\address{Institut Gaspard Monge
UMR CNRS - 8049
Universit\'e Gustave Eiffel
5 boulevard Descartes, 
77454 Champs-Sur-Marne
FRANCE}

\begin{abstract}
We give a general construction of triangulations starting from a walk in the quarter plane with small steps, which is a discrete version of the mating of trees. We use a special instance of this construction to give a bijection between maps equipped with a rooted spanning tree and walks in the quarter plane.
We also show how the construction allows to
  recover several known bijections between such objects  in a uniform way.
 
\end{abstract}
\maketitle
\section{Introduction}

Mating of polynomials originates in complex dynamics, where one can match two Julia sets in order to build a topological sphere or a surface,
see e.g.  \cite{BEKMPRT} and \newline
\centerline{\tt https://www.math.univ-toulouse.fr/~cheritat/MatMovies/}\newline
for nice pictures and movies. This includes in particular the case of Julia sets which are topologically real trees.
This construction  has been introduced in probability by Le Gall and Paulin \cite{LP} for studying the topology of the Brownian map and then used,
under the name ``mating of trees'' by 
 Duplantier, Miller  and Sheffield \cite{DMS} in quantum gravity, followed by many others.
A discrete version of this construction already appears in Mullin's bijection (see \cite{Sch97}) and has been further used recently by several authors, see e.g.  Gwynne, Holden and  Sun \cite{GHS} for a recent overview.
The basic idea is to consider a walk in the quarter plane $x,y\geq 0$, with steps in the set $\{(1,0),(-1,0),(0,1),(0,-1)\}$, starting and ending at $0$, and to write the vertical and horizontal coordinates of the walk as two opposite Motzkin paths running vertically. See   below, with one Motzkin path on the left for the horizontal coordinate and one on the right for the vertical coordinate:

$${\tt    \setlength{\unitlength}{1pt}
\begin{picture}(100,130)
\thicklines    
\

\put(0,0){\line(1,0){80}}
\put(0,0){\line(0,1){120}}
\put(0,120){\line(1,0){80}}
\put(80,0){\line(0,1){120}}

\put(0,0){\color{red}\line(0,1){10}}             \put(80,0){\color{blue}\line(-1,1){10}}
\put(0,10){\color{red}\line(1,1){10}}            \put(70,10){\color{blue}\line(0,1){10}}
\put(10,20){\color{red}\line(0,1){10}}          \put(70,20){\color{blue}\line(-1,1){10}} 
\put(10,30){\color{red}\line(-1,1){10}}          \put(60,30){\color{blue}\line(0,1){20}} 
\put(0,40){\color{red}\line(1,1){10}}           \put(60,50){\color{blue}\line(1,1){10}}
\put(10,50){\color{red}\line(0,1){10}}            \put(70,60){\color{blue}\line(0,1){10}}    
\put(10,60){\color{red}\line(1,1){10}}        \put(70,70){\color{blue}\line(1,1){10}} 
\put(20,70){\color{red}\line(0,1){20}}          \put(80,80){\color{blue}\line(-1,1){10}} 
\put(20,90){\color{red}\line(-1,1){20}}           \put(70,90){\color{blue}\line(0,1){20}}
\put(0,110){\color{red}\line(0,1){10}}            \put(70,110){\color{blue}\line(1,1){10}}    

\end{picture}}
$$

then mate the two Motzkin paths by drawing horizontal lines between vertices:
$${\tt    \setlength{\unitlength}{1pt}
\begin{picture}(100,130)
\thicklines    
\

\put(0,0){\line(1,0){80}}
\put(0,0){\line(0,1){120}}
\put(0,120){\line(1,0){80}}
\put(80,0){\line(0,1){120}}

\put(0,0){\color{red}\line(0,1){10}}             \put(80,0){\color{blue}\line(-1,1){10}}
\put(0,10){\color{red}\line(1,1){10}}            \put(70,10){\color{blue}\line(0,1){10}}
\put(10,20){\color{red}\line(0,1){10}}          \put(70,20){\color{blue}\line(-1,1){10}} 
\put(10,30){\color{red}\line(-1,1){10}}          \put(60,30){\color{blue}\line(0,1){20}} 
\put(0,40){\color{red}\line(1,1){10}}           \put(60,50){\color{blue}\line(1,1){10}}
\put(10,50){\color{red}\line(0,1){10}}            \put(70,60){\color{blue}\line(0,1){10}}    
\put(10,60){\color{red}\line(1,1){10}}        \put(70,70){\color{blue}\line(1,1){10}} 
\put(20,70){\color{red}\line(0,1){20}}          \put(80,80){\color{blue}\line(-1,1){10}} 
\put(20,90){\color{red}\line(-1,1){20}}           \put(70,90){\color{blue}\line(0,1){20}}
\put(0,110){\color{red}\line(0,1){10}}            \put(70,110){\color{blue}\line(1,1){10}}    


\put(0,10){\line(1,0){70}}  
\put(10,20){\line(1,0){60}}  
\put(10,30){\line(1,0){50}} 
\put(0,40){\line(1,0){60}} 
\put(10,50){\line(1,0){50}} 
\put(10,60){\line(1,0){60}} 
\put(20,70){\line(1,0){50}} 
\put(20,80){\line(1,0){60}} 
\put(20,90){\line(1,0){50}} 
\put(10,100){\line(1,0){60}} 
\put(0,110){\line(1,0){70}}

\end{picture}}
$$
The Motzkin paths are then contracted into trees in the usual way, while the upper and lower boundaries of the rectangle are identified.
The resulting map is a planar triangulation. I will give a more detailed explanation in the following sections.

In this note I use a generalization of these ideas to walks with {\sl small steps}, i.e. with steps taken in the set 
$\{(1,0),(-1,0),(0,1),(0,-1),(-1,-1),(1,-1),(-1,1),(1,1)\}$, which involves constructing a map with faces of degrees three and four then contracting the faces of degree four. As we shall see, the generality of the construction and the many variants one can produce from  it make it a versatile tool for producing bijections between specific classes of walks 
and of maps. Again, the precise definitions will be given in the next section.
Using this I give a bijection between a certain class of walks in the quarter plane (which I call reversed $Y$-walks, see below) and maps equipped with a rooted spanning tree in which the degrees of the vertices are encoded by the length of the steps of the walk.
This bijection is quite different from Mullin's bijection but bears some connection with blossoming bijections \cite{Sch97}. These ideas also
allow us to reinterpret several known bijections between walks and triangulations, or more general planar maps. In particular we will show how  the following bijections fit into this framework:

\begin{itemize}
\item
The bijection
of Mullin, between walks in the quarter plane with straight steps and triangulations having a Hamiltonian cycle.
\item A bijection of Bernardi between Kreweras walks and triangulations with a spanning tree.
\item
A bijection  of Kenyon, Miller, Sheffield and Wilson between  walks and maps with a bipolar orientation.
\item A bijection 
 of Li, Sun and Watson between tandem walks satisfying some further conditions and Schnyder woods.
\end{itemize}
Rather than the particular results which are obtained, we think that the main interest of this paper is its general philosophy and its potential to produce a wealth of bijections between walks and maps.

This paper is organized as follows. In the next section I give a general algorithm for producing  planar  triangulations, starting from a walk with small steps in the quarter plane. In order to obtain precise bijections one needs to specify a number of parameters in this construction.  In section \ref{sec:span} I introduce reversed $Y$-walks in the quarter plane and I give a new bijection between these  walks and planar maps equipped with a spanning tree. Then, in section \ref{sec:4}, I show how to recover the bijections listed above using similar considerations.

\section{Associating a triangulation  to a walk in the quarter plane}
\label{sec:3}

\subsection{Trees and Dyck paths}\label{sec:2}

There is a well known way to  associate a rooted planar tree to a Dyck path by 
  matching up and down steps:

$${\tt    \setlength{\unitlength}{1.5pt}
\begin{picture}(200,40)
\thicklines    

\put(0,-10){}
\put(0,0){\line(1,0){100}}
\put(0,0){\color{blue}\line(1,1){10}}
\put(10,10){\color{blue}\line(1,1){10}}
\put(20,20){\color{blue}\line(1,-1){10}}
\put(30,10){\color{blue}\line(1,1){10}}
\put(40,20){\color{blue}\line(1,1){10}}
\put(50,30){\color{blue}\line(1,-1){10}}
\put(60,20){\color{blue}\line(1,1){10}}
\put(70,30){\color{blue}\line(1,-1){10}}
\put(80,20){\color{blue}\line(1,-1){10}}
\put(90,10){\color{blue}\line(1,-1){10}}
\put(115,20){$\to$}
\put(148,-2){\textcolor{blue}{$\bullet$}}
\put(150,0){\color{blue}\line(1,1){10}}
\put(158,8){\textcolor{blue}{$\bullet$}}
\put(160,10){\color{blue}\line(1,1){10}}
\put(160,10){\color{blue}\line(-1,1){10}}
\put(148,18){\textcolor{blue}{$\bullet$}}
\put(168,18){\textcolor{blue}{$\bullet$}}
\put(170,20){\color{blue}\line(1,1){10}}
\put(170,20){\color{blue}\line(-1,1){10}}
\put(158,28){\textcolor{blue}{$\bullet$}}
\put(178,28){\textcolor{blue}{$\bullet$}}

\thinlines    
\put(2.5,2.5){\line(1,0){95}}\put(5,5){\line(1,0){90}}\put(7.5,7.5){\line(1,0){85}}\put(10,10){\line(1,0){80}}
\put(12.5,12.5){\line(1,0){15}}\put(15,15){\line(1,0){10}}\put(17.5,17.5){\line(1,0){5}}
\put(32.5,12.5){\line(1,0){55}}\put(35,15){\line(1,0){50}}\put(37.5,17.5){\line(1,0){45}}\put(40,20){\line(1,0){40}}
\put(42.5,22.5){\line(1,0){15}}\put(45,25){\line(1,0){10}}\put(47.5,27.5){\line(1,0){5}}
\put(62.5,22.5){\line(1,0){15}}\put(65,25){\line(1,0){10}}\put(67.5,27.5){\line(1,0){5}}
\end{picture}}
$$

One can think that one cuts out from the plane the striped area under the Dyck path and sews the up and down steps to produce the tree embedded into the plane.
One can generalize this  construction to Motzkin paths by shrinking each horizontal step to a point.

$${\tt    \setlength{\unitlength}{1.5pt}
\begin{picture}(200,40)
\thicklines    

\put(0,-10){}
\put(0,0){\line(1,0){110}}
\put(0,0){\color{blue}\line(1,1){10}}
\put(10,10){\color{blue}\line(1,1){10}}
\put(20,20){\color{blue}\line(1,-1){10}}
\put(30,10){\color{blue}\line(1,1){10}}
\put(40,20){\color{blue}\line(1,1){10}}
\put(50,30){\color{blue}\line(1,-1){10}}
\put(60,20){\color{blue}\line(1,0){10}}
\put(70,20){\color{blue}\line(1,1){10}}
\put(80,30){\color{blue}\line(1,-1){10}}
\put(90,20){\color{blue}\line(1,-1){10}}
\put(100,10){\color{blue}\line(1,-1){10}}
\put(115,20){$\to$}
\put(148,-2){\textcolor{blue}{$\bullet$}}
\put(150,0){\color{blue}\line(1,1){10}}
\put(158,8){\textcolor{blue}{$\bullet$}}
\put(160,10){\color{blue}\line(1,1){10}}
\put(160,10){\color{blue}\line(-1,1){10}}
\put(148,18){\textcolor{blue}{$\bullet$}}
\put(168,18){\textcolor{blue}{$\bullet$}}
\put(170,20){\color{blue}\line(1,1){10}}
\put(170,20){\color{blue}\line(-1,1){10}}
\put(158,28){\textcolor{blue}{$\bullet$}}
\put(178,28){\textcolor{blue}{$\bullet$}}

\thinlines    
\put(2.5,2.5){\line(1,0){105}}\put(5,5){\line(1,0){100}}\put(7.5,7.5){\line(1,0){95}}\put(10,10){\line(1,0){90}}
\put(12.5,12.5){\line(1,0){15}}\put(15,15){\line(1,0){10}}\put(17.5,17.5){\line(1,0){5}}
\put(32.5,12.5){\line(1,0){65}}\put(35,15){\line(1,0){60}}\put(37.5,17.5){\line(1,0){55}}\put(40,20){\line(1,0){50}}
\put(42.5,22.5){\line(1,0){15}}\put(45,25){\line(1,0){10}}\put(47.5,27.5){\line(1,0){5}}
\put(72.5,22.5){\line(1,0){15}}\put(75,25){\line(1,0){10}}\put(77.5,27.5){\line(1,0){5}}
\end{picture}}
$$

\subsection{}\label{sec2:2}
One can also consider paths which do not start or end at $0$
and build a forest of trees, rooted on a $V$-shaped path:

$${\tt    \setlength{\unitlength}{1.5pt}
\begin{picture}(250,40)
\thicklines    

\put(0,-10){}
\put(0,0){\line(1,0){130}}
\put(0,20){\color{blue}\line(1,1){10}}
\put(10,30){\color{blue}\line(1,-1){10}}
\put(20,20){\color{blue}\line(1,1){10}}
\put(30,30){\color{blue}\line(1,-1){10}}
\put(40,20){\color{red}\line(1,-1){10}}
\put(50,10){\color{red}\line(1,-1){10}}
\put(60,0){\color{blue}\line(1,1){10}}
\put(70,10){\color{blue}\line(1,-1){10}}
\put(80,0){\color{red}\line(1,1){10}}
\put(90,10){\color{red}\line(1,1){10}}
\put(100,20){\color{blue}\line(1,1){10}}
\put(110,30){\color{blue}\line(1,-1){10}}
\put(120,20){\color{red}\line(1,1){10}}
\put(150,20){$\to$}
\put(198,-2){\textcolor{blue}{$\bullet$}}
\put(200,0){\color{red}\line(1,1){30}}
\put(170,30){\color{blue}\line(1,-1){10}}
\put(180,20){\color{red}\line(1,-1){20}}
\put(178,18){\textcolor{blue}{$\bullet$}}
\put(198,8){\textcolor{blue}{$\bullet$}}
\put(200,0){\color{blue}\line(0,1){10}}
\put(168,28){\textcolor{blue}{$\bullet$}}
\put(188,28){\textcolor{blue}{$\bullet$}}
\put(180,20){\color{blue}\line(1,1){10}}
\put(220,20){\color{blue}\line(-1,1){10}}
\put(218,18){\textcolor{blue}{$\bullet$}}
\put(208,28){\textcolor{blue}{$\bullet$}}

\thinlines    
\put(0,20){\line(1,0){40}}
\put(2.5,22.5){\line(1,0){15}}\put(5,25){\line(1,0){10}}\put(7.5,27.5){\line(1,0){5}}
\put(22.5,22.5){\line(1,0){15}}\put(25,25){\line(1,0){10}}\put(27.5,27.5){\line(1,0){5}}

\put(62.5,2.5){\line(1,0){15}}\put(65,5){\line(1,0){10}}\put(67.5,7.5){\line(1,0){5}}
\put(102.5,22.5){\line(1,0){15}}\put(105,25){\line(1,0){10}}\put(107.5,27.5){\line(1,0){5}}
\put(100,20){\line(1,0){20}}

\put(180,0){\line(1,0){50}}

\end{picture}}
$$

Some down steps (before the minimum of the path is reached) or upward steps (after the minimum) are left unmatched and are indicated in red on the picture. They form the ``$V$'' on which the forest is rooted.

\subsection{The basic construction}\label{sec:3.1}
We consider walks in the quarter plane with   {\sl small steps}:  the  coordinates of the steps belong to  $\{-1,0,1\}$, thus we get $8$ distinct non zero steps:

$${\tt    \setlength{\unitlength}{1.5pt}
\begin{picture}(30,20)
\thicklines    

\put(0,10){\vector(1,0){10}}
\put(0,10){\vector(1,1){10}}
\put(0,10){\vector(-1,0){10}}
\put(0,10){\vector(1,-1){10}}
\put(0,10){\vector(-1,1){10}}
\put(0,10){\vector(-1,-1){10}}
\put(0,10){\vector(0,1){10}}
\put(0,10){\vector(0,-1){10}}
\end{picture}}
$$
We can distinguish straight steps:
$(1,0),(0,1),(-1,0),(0,-1)$ 
$${\tt    \setlength{\unitlength}{1.5pt}
\begin{picture}(30,20)
\thicklines

\put(0,10){\vector(1,0){10}}

\put(0,10){\vector(-1,0){10}}

\put(0,10){\vector(0,1){10}}
\put(0,10){\vector(0,-1){10}}
\end{picture}}
$$
from oblique steps:
$${\tt    \setlength{\unitlength}{1.5pt}
\begin{picture}(30,20)
\thicklines

\put(0,10){\vector(1,1){10}}

\put(0,10){\vector(1,-1){10}}
\put(0,10){\vector(-1,1){10}}
\put(0,10){\vector(-1,-1){10}}

\end{picture}}
$$ 

The problem of enumerating such walks has been considered in many papers, starting from the seminal work \cite{BMM}.

Consider a two-dimensional walk with small steps, starting and ending at  the origin, which remains in the quarter plane $x,y\geq 0$. 
Here is an example with twelve steps, ten of which are straight and two oblique:
$$(1,0);(0,1);(-1,1);(1,0);(0,-1);(1,0);(0,-1);(0,1);(-1,0);(-1,0);(1,0);(-1,-1).$$
The walk is shown below:

$$
\includegraphics[width=4cm]{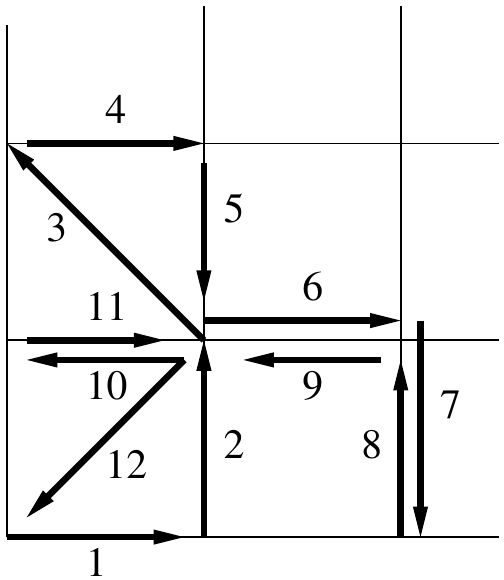}
  $$

The projections of this walk on the coordinate axes are two Motzkin paths with respective~steps
$1;0;-1;1;0;1;0;0;-1;-1;1;-1$ and $0;1;1;0;-1;0;-1;1;0;0;0;-1$:
$${\tt    \setlength{\unitlength}{1pt}
\begin{picture}(300,40)
\thicklines    
\

\put(0,0){\line(1,0){120}}
\put(150,0){\line(1,0){120}}

\put(0,0){\color{red}\line(1,1){10}}            \put(150,0){\color{blue}\line(1,0){10}}
\put(10,10){\color{red}\line(1,0){10}}          \put(160,0){\color{blue}\line(1,1){20}} 
\put(20,10){\color{red}\line(1,-1){10}}          \put(180,20){\color{blue}\line(1,0){10}} 
\put(30,0){\color{red}\line(1,1){10}}           \put(190,20){\color{blue}\line(1,-1){10}}
\put(40,10){\color{red}\line(1,0){10}}            \put(200,10){\color{blue}\line(1,0){10}}    
\put(50,10){\color{red}\line(1,1){10}}        \put(210,10){\color{blue}\line(1,-1){10}} 
\put(60,20){\color{red}\line(1,0){20}}          \put(220,0){\color{blue}\line(1,1){10}} 
\put(80,20){\color{red}\line(1,-1){20}}           \put(230,10){\color{blue}\line(1,0){30}}
\put(100,0){\color{red}\line(1,1){10}}            \put(260,10){\color{blue}\line(1,-1){10}}    
\put(110,10){\color{red}\line(1,-1){10}}        
\end{picture}}
$$

We draw the two paths vertically in an  opposite way,
with the horizontal coordinate on the left and the vertical coordinate on the right, the paths running from bottom to top.
 
$$\tt    \setlength{\unitlength}{1pt}
\begin{picture}(100,130)
\thicklines    
\

\put(0,10){\line(1,0){80}}
\put(0,10){\line(0,1){120}}
\put(0,130){\line(1,0){80}}
\put(80,10){\line(0,1){120}}

\put(0,10){\color{red}\line(1,1){10}}            \put(80,10){\color{blue}\line(0,1){10}}
\put(10,20){\color{red}\line(0,1){10}}          \put(80,20){\color{blue}\line(-1,1){20}} 
\put(10,30){\color{red}\line(-1,1){10}}          \put(60,40){\color{blue}\line(0,1){10}} 
\put(0,40){\color{red}\line(1,1){10}}           \put(60,50){\color{blue}\line(1,1){10}}
\put(10,50){\color{red}\line(0,1){10}}            \put(70,60){\color{blue}\line(0,1){10}}    
\put(10,60){\color{red}\line(1,1){10}}        \put(70,70){\color{blue}\line(1,1){10}} 
\put(20,70){\color{red}\line(0,1){20}}          \put(80,80){\color{blue}\line(-1,1){10}} 
\put(20,90){\color{red}\line(-1,1){20}}           \put(70,90){\color{blue}\line(0,1){30}}
\put(0,110){\color{red}\line(1,1){10}}            \put(70,120){\color{blue}\line(1,1){10}}    
\put(10,120){\color{red}\line(-1,1){10}}        
\end{picture}
$$

The mating consists in drawing horizontal lines between the vertices of the two paths, as below:

$${\tt    \setlength{\unitlength}{1pt}
\begin{picture}(100,130)
\thicklines    
\

\put(0,10){\line(1,0){80}}
\put(0,10){\line(0,1){120}}
\put(0,130){\line(1,0){80}}
\put(80,10){\line(0,1){120}}

\put(0,10){\color{red}\line(1,1){10}}            \put(80,10){\color{blue}\line(0,1){10}}
\put(10,20){\color{red}\line(0,1){10}}          \put(80,20){\color{blue}\line(-1,1){20}} 
\put(10,30){\color{red}\line(-1,1){10}}          \put(60,40){\color{blue}\line(0,1){10}} 
\put(0,40){\color{red}\line(1,1){10}}           \put(60,50){\color{blue}\line(1,1){10}}
\put(10,50){\color{red}\line(0,1){10}}            \put(70,60){\color{blue}\line(0,1){10}}    
\put(10,60){\color{red}\line(1,1){10}}        \put(70,70){\color{blue}\line(1,1){10}} 
\put(20,70){\color{red}\line(0,1){20}}          \put(80,80){\color{blue}\line(-1,1){10}} 
\put(20,90){\color{red}\line(-1,1){20}}           \put(70,90){\color{blue}\line(0,1){30}}
\put(0,110){\color{red}\line(1,1){10}}            \put(70,120){\color{blue}\line(1,1){10}}    
\put(10,120){\color{red}\line(-1,1){10}}

\put(10,10){\line(1,0){70}} 
\put(10,20){\line(1,0){70}}
\put(10,30){\line(1,0){60}} 
\put(0,40){\line(1,0){60}} 
\put(10,50){\line(1,0){50}} 
\put(10,60){\line(1,0){60}} 
\put(20,70){\line(1,0){50}} 
 \put(20,80){\line(1,0){60}} 
\put(20,90){\line(1,0){50}} 
\put(10,100){\line(1,0){60}} 
\put(0,110){\line(1,0){70}} 
 \put(10,120){\line(1,0){60}} 
 
\end{picture}}
$$

Between the two Motzkin paths we have then a succession of quadrilaterals, each one corresponding to a step of the walk. These quadrilaterals have several types,  
some of them, corresponding to straight steps,  have one of their side vertical, while the others, corresponding to oblique steps, have tilted sides. Here are the straight steps, numbered from bottom to top, and the shaded oblique steps.

$${\tt    \setlength{\unitlength}{1pt}
\begin{picture}(100,130)
\thicklines    
\

\put(0,10){\line(1,0){80}}
\put(0,10){\line(0,1){120}}
\put(0,130){\line(1,0){80}}
\put(80,10){\line(0,1){120}}

\put(0,10){\color{red}\line(1,1){10}}            \put(80,10){\color{blue}\line(0,1){10}}
\put(10,20){\color{red}\line(0,1){10}}          \put(80,20){\color{blue}\line(-1,1){20}} 
\put(10,30){\color{red}\line(-1,1){10}}          \put(60,40){\color{blue}\line(0,1){10}} 
\put(0,40){\color{red}\line(1,1){10}}           \put(60,50){\color{blue}\line(1,1){10}}
\put(10,50){\color{red}\line(0,1){10}}            \put(70,60){\color{blue}\line(0,1){10}}    
\put(10,60){\color{red}\line(1,1){10}}        \put(70,70){\color{blue}\line(1,1){10}} 
\put(20,70){\color{red}\line(0,1){20}}          \put(80,80){\color{blue}\line(-1,1){10}} 
\put(20,90){\color{red}\line(-1,1){20}}           \put(70,90){\color{blue}\line(0,1){30}}
\put(0,110){\color{red}\line(1,1){10}}            \put(70,120){\color{blue}\line(1,1){10}}    
\put(10,120){\color{red}\line(-1,1){10}}     

\put(36,12){\text{$\scriptstyle 1$}}
\put(36,22){\text{$\scriptstyle 2$}}
\put(36,42){\text{$\scriptstyle 3$}}
\put(36,52){\text{$\scriptstyle 4$}}
\put(36,62){\text{$\scriptstyle 5$}}
\put(36,72){\text{$\scriptstyle 6$}}
\put(36,82){\text{$\scriptstyle 7$}}
\put(36,92){\text{$\scriptstyle 8$}}
\put(36,102){\text{$\scriptstyle 9$}}
\put(34,112){\text{$\scriptstyle 10$}}

\put(10,10){\line(1,0){70}} 
\put(10,20){\line(1,0){70}}
\put(10,30){\line(1,0){60}} 
\put(0,40){\line(1,0){60}} 
\put(10,50){\line(1,0){50}} 
\put(10,60){\line(1,0){60}} 
\put(20,70){\line(1,0){50}} 
 \put(20,80){\line(1,0){60}} 
\put(20,90){\line(1,0){50}} 
\put(10,100){\line(1,0){60}} 
\put(0,110){\line(1,0){70}} 
 \put(10,120){\line(1,0){60}} 

\put(8,32){\color{green}\line(1,0){60}}        
\put(6,34){\color{green}\line(1,0){60}}
\put(4,36){\color{green}\line(1,0){60}}
\put(2,38){\color{green}\line(1,0){60}}
\put(8,122){\color{green}\line(1,0){64}}        
\put(6,124){\color{green}\line(1,0){68}}
\put(4,126){\color{green}\line(1,0){72}}
\put(2,128){\color{green}\line(1,0){76}}
\end{picture}}
$$

We now contract the two Motzkin paths into two trees as in section \ref{sec:2}. Again we can imagine that we cut out the area below each Motzkin path and sew the boundaries.
In this contraction  the quadrilaterals corresponding to straight steps become triangles as, for example:

$${\tt    \setlength{\unitlength}{1pt}
\begin{picture}(300,15)
\thinlines    
\

\put(0,0){\line(1,0){60}}
\put(60,0){\color{blue}\line(0,1){10}}
\put(10,10){\line(1,0){50}}

\put(0,0){\color{red}\line(1,1){10}}           
\put (80,5){$\to$}
\put(100,0){\line(1,0){60}}

\put(110,10){\line(5,-1){50}}

\put(100,0){\color{red}\line(1,1){10}}           

\end{picture}}
$$

\begin{figure}
  
  \caption{Contraction of the Motzkin paths.}
  \label{fig:contraction1}

$$
\includegraphics[width=8cm]{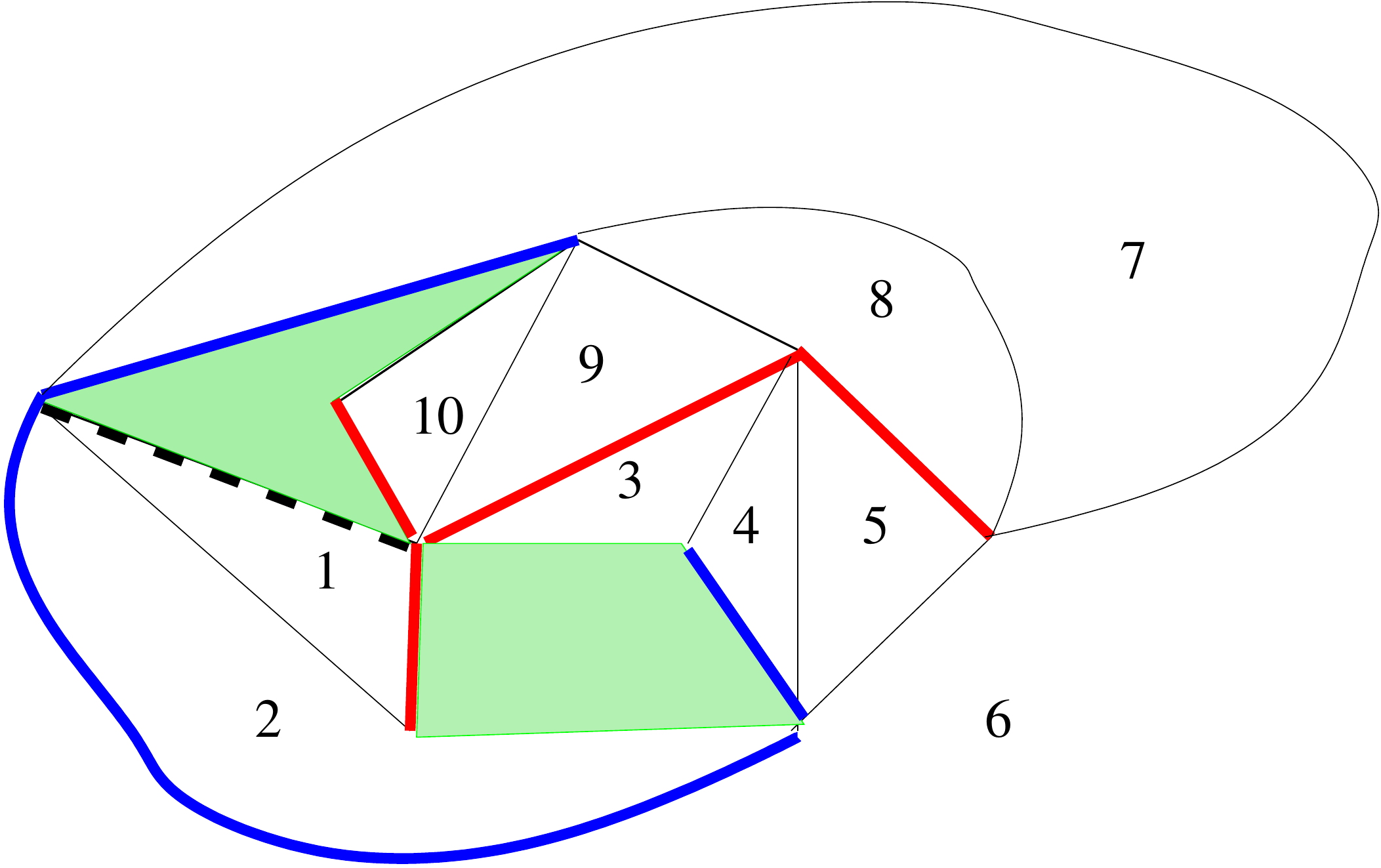}
  $$
\end{figure}

Once we have made these contractions we obtain 
a planar map  whose faces are triangles (corresponding to straight steps), quadrilaterals (corresponding to oblique steps) and an external face with two sides, the two horizontal sides of the initial  rectangle.
If we identify these two sides by contracting the external face, we  obtain a planar map, with a distinguished edge, the one corresponding to the two sides of the rectangle,   on which two trees are drawn, namely the images of the two Motzkin paths.
Figure \ref{fig:contraction1} shows the map obtained from the above walk, where the triangles are numbered  and the two quadrilaterals are shaded, while the two trees are depicted in red and blue and the distinguished edge (corresponding to the upper and lower boundaries of the rectangle) is dashed.

In order to obtain a planar triangulation we will contract the quadrilaterals corresponding to oblique steps: again imagine that one cuts out one of these quadrilaterals from the plane, then   identifies two opposite vertices, pairs the edges correspondingly and sews them. 
More precisely, consider such a quadrilateral:

$${\tt    \setlength{\unitlength}{1pt}
\begin{picture}(300,15)
\thinlines    
\

\put(0,0){\line(1,0){60}}
\put(60,0){\color{blue}\line(1,1){10}}
\put(-10,10){\line(1,0){80}}

\put(0,0){\color{red}\line(-1,1){10}}           

\end{picture}}
$$
There are two ways to contract  this quadrilateral along one of its diagonals to produce two segments, either like this by identifying opposite vertices 
 $v$ and $w$ and accordingly sewing the edge $uv$ with $uw$ and the edge $vz$ with $wz$:
$${\tt    \setlength{\unitlength}{1pt}
\begin{picture}(300,30)
\thinlines    

\put(-8,-8){$w$}\put(63,-8){$z$}\put(-16,16){$u$}\put(73,16){$v$}
\put(0,0){\line(1,0){60}}
\put(60,0){\color{blue}\line(1,1){10}}
\put(-10,10){\line(1,0){80}}

\put(0,0){\color{red}\line(-1,1){10}}           
 \put(80,5){$\to$}
\put(100,5){$u$}
\put(110,5){\color{red}\line(1,0){40}}\put(150,5){\color{blue}\line(1,0){40}}\put(148,8){$v$}\put(148,-4){$w$}\put(193,5){$z$}
\end{picture}}
$$

\medskip

or like this, identifying  $u$ and $z$:
$${\tt    \setlength{\unitlength}{1pt}
\begin{picture}(300,30)
\thinlines    
\
\put(-6,-8){$w$}\put(63,-8){$z$}\put(-16,16){$u$}\put(73,16){$v$}
\put(0,0){\line(1,0){60}}
\put(60,0){\color{blue}\line(1,1){10}}
\put(-10,10){\line(1,0){80}}

\put(0,0){\color{red}\line(-1,1){10}}           
 \put(80,5){$\to$}
\put(100,5){$w$}
\put(110,5){\color{red}\line(1,0){40}}\put(150,5){\color{blue}\line(1,0){40}}\put(148,8){$u$}\put(148,-4){$z$}\put(193,5){$v$}
\end{picture}}
$$

\medskip

%

\begin{figure}
  
  \caption{Contracting the quadrilaterals of Figure \ref{fig:contraction1}.}
  \label{fig:contraction2}

$$
\includegraphics[width=8cm]{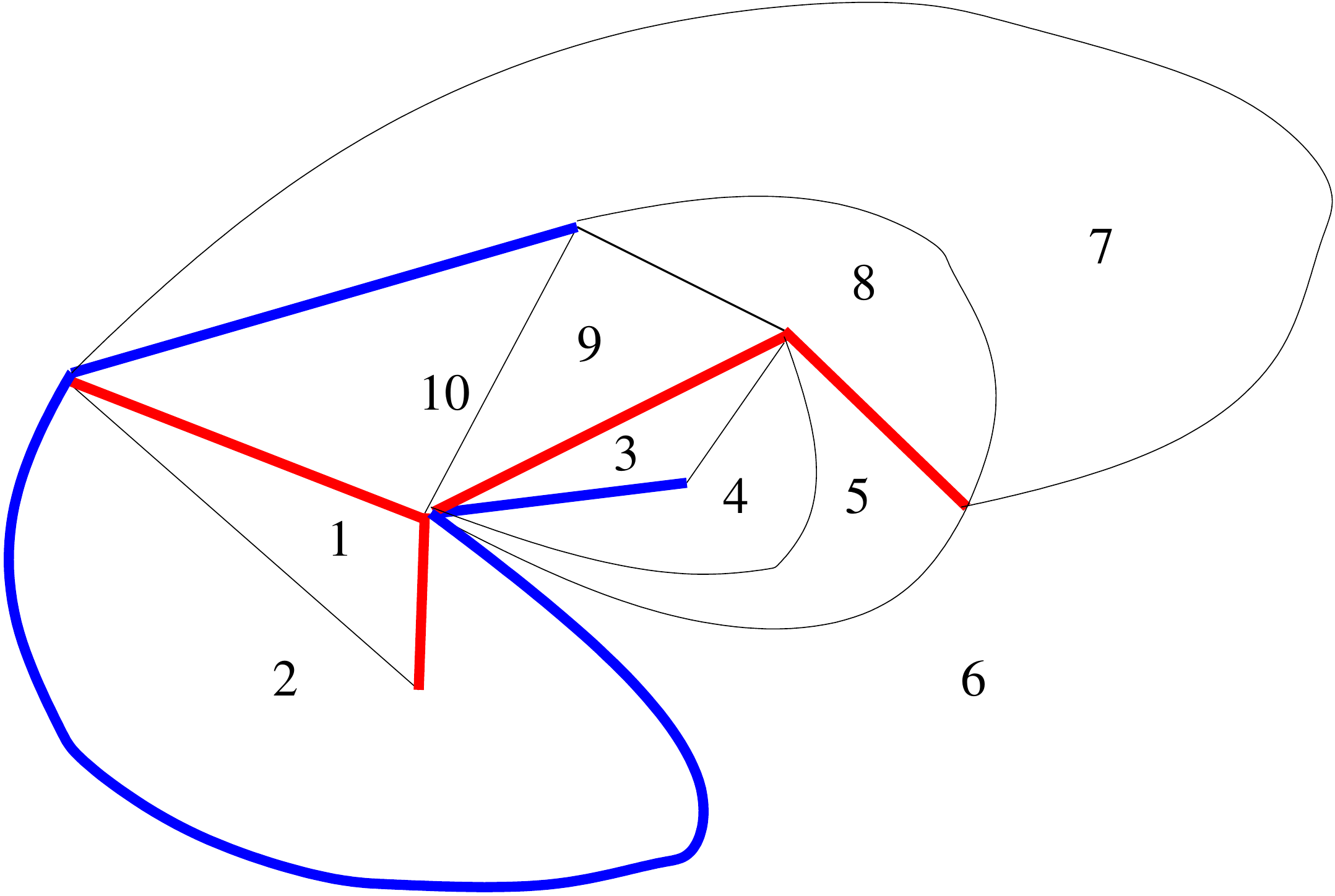}
  $$
\end{figure}

Let $o$ be the number of oblique steps in the walk.
From the map corresponding to the walk we can obtain  $2^{o}$ triangulations by removing the  quadrilaterals and, for each one, sewing its boundaries according to the two possibilites. Figure \ref{fig:contraction2}  shows the result of a choice of such contractions on the map of Figure \ref{fig:contraction1}.

 The construction presented above is very general, but  not one-to-one. Indeed one can associate to each path $2^{o}$ triangulations,
 moreover some triangulations may be obtained by more than one of these constructions.

\begin{figure}
$$
\includegraphics[width=5cm]{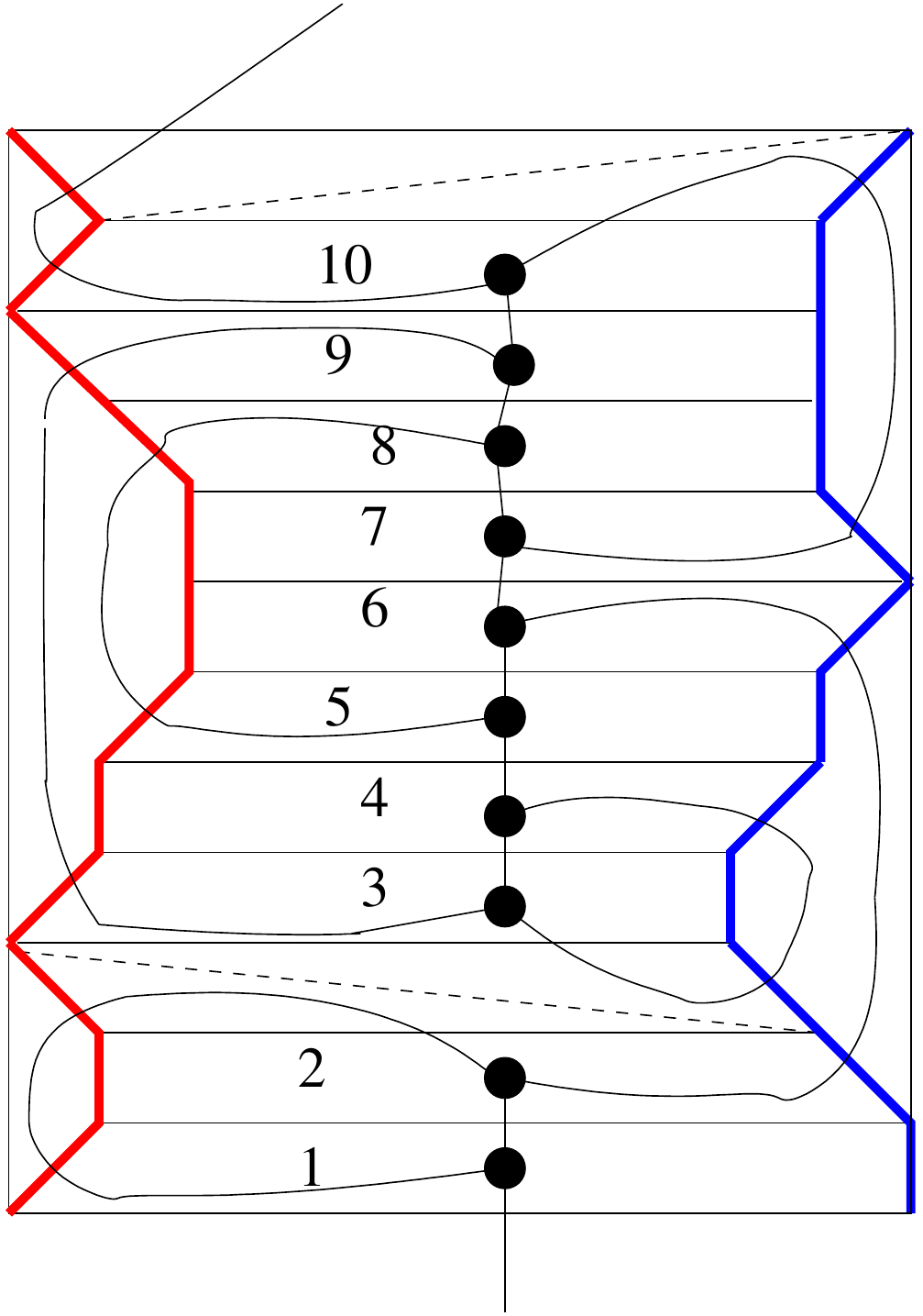}
  $$
\caption{The triangulation and its dual cubic map}
\label{fig:dual}
\end{figure}

 In the following I will consider specific instances of the  construction, where one restricts the class of steps which are available for the walk and one gives an explicit algorithm to decide, for each quadrilateral,  which of the two contractions is made. 
In order to visualize the maps obtained by these constructions we will find it useful to depict also the dual map as in Figure \ref{fig:dual},
where we indicate the way we contract the quadrilateral by drawing a dashed line between the opposite vertices which have been identified, then draw edges connecting the dual vertices.

Before going into the description of the examples,  I describe some possible variations on this construction.
\subsection{Variations}\label{sec:3:2}
\subsubsection{Changing the end point}
It is possible to generalize the construction to paths which do not start or end at zero.  Also one can relax the condition that the walk remains in the quarter plane. In this case the two paths are contracted as explained in  section \ref{sec2:2} and there remain edges which are not matched. Together with the upper and lower horizontal sides, they form an external boundary. Here is an example, where we show vertical lines between matched edges of the paths. There remains three unmatched edges on the left and two on the right:
$$
\includegraphics[width=4cm]{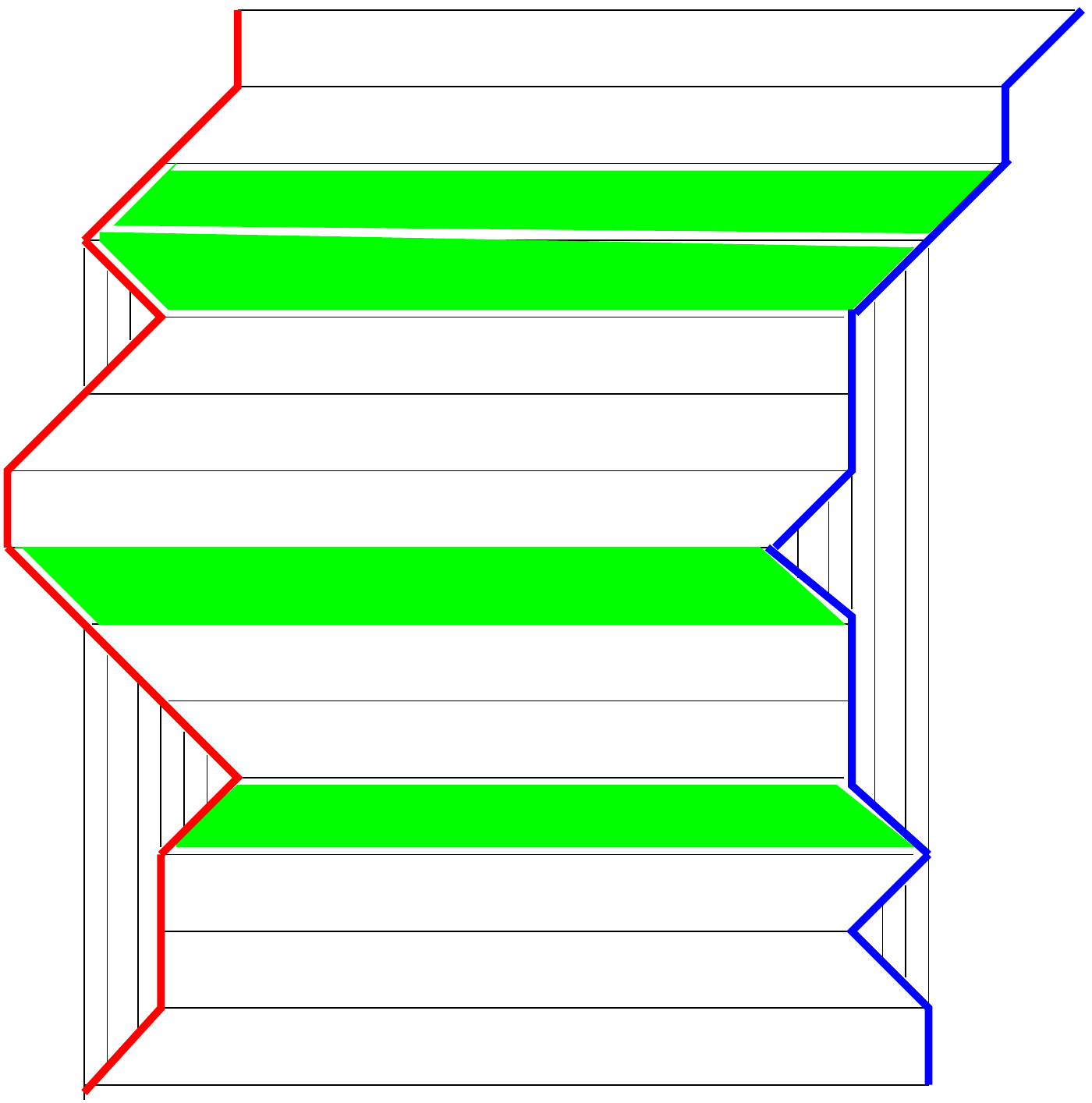}
  $$

At this point, one can either keep the unmatched edges and consider that they bound an exterior face, or  idenfity  
 these edges pairwise, if there is an even number of them. 
 We will see some examples in the following.

\subsubsection{Changing the set of steps}\label{step}
It is also possible to use other types of steps. If one uses steps of type $(i,j)$ then one can consider the path obtained by replacing the step $(i,j)$ by a sequence of $|i|$  steps of type $(\sgn(i),0)$ and $|j|$ steps of type $(0,\sgn(j))$, then erase the horizontal lines in the polygon thus formed to get a face with $|i|+|j|+2$ edges on the boundary. For example here, in the case $i=3,j=4$, we get the following picture:
$$
\includegraphics[width=5cm]{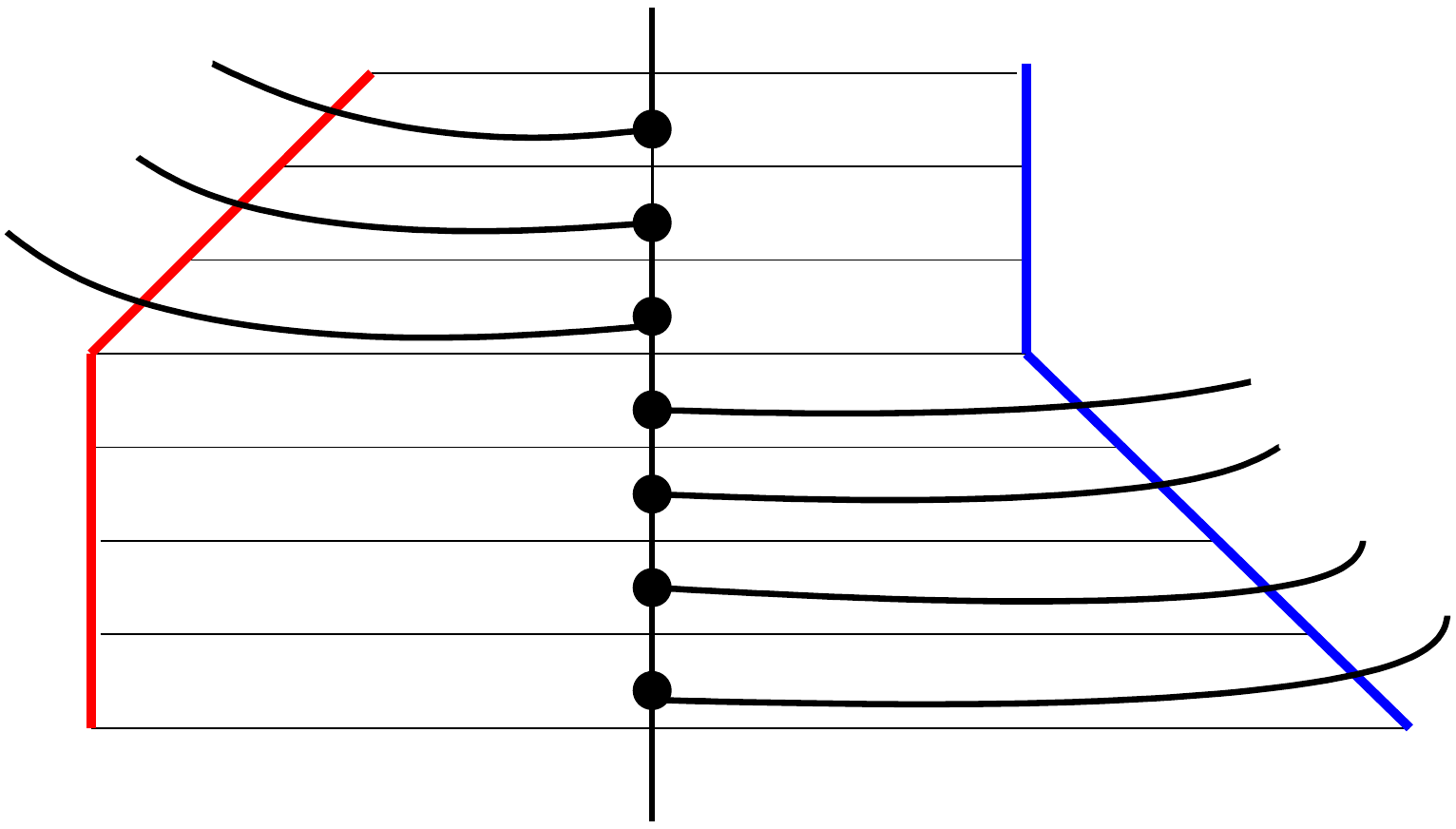}
  $$
When we erase the internal edges we get  a face with nine edges
 $$
\includegraphics[width=5cm]{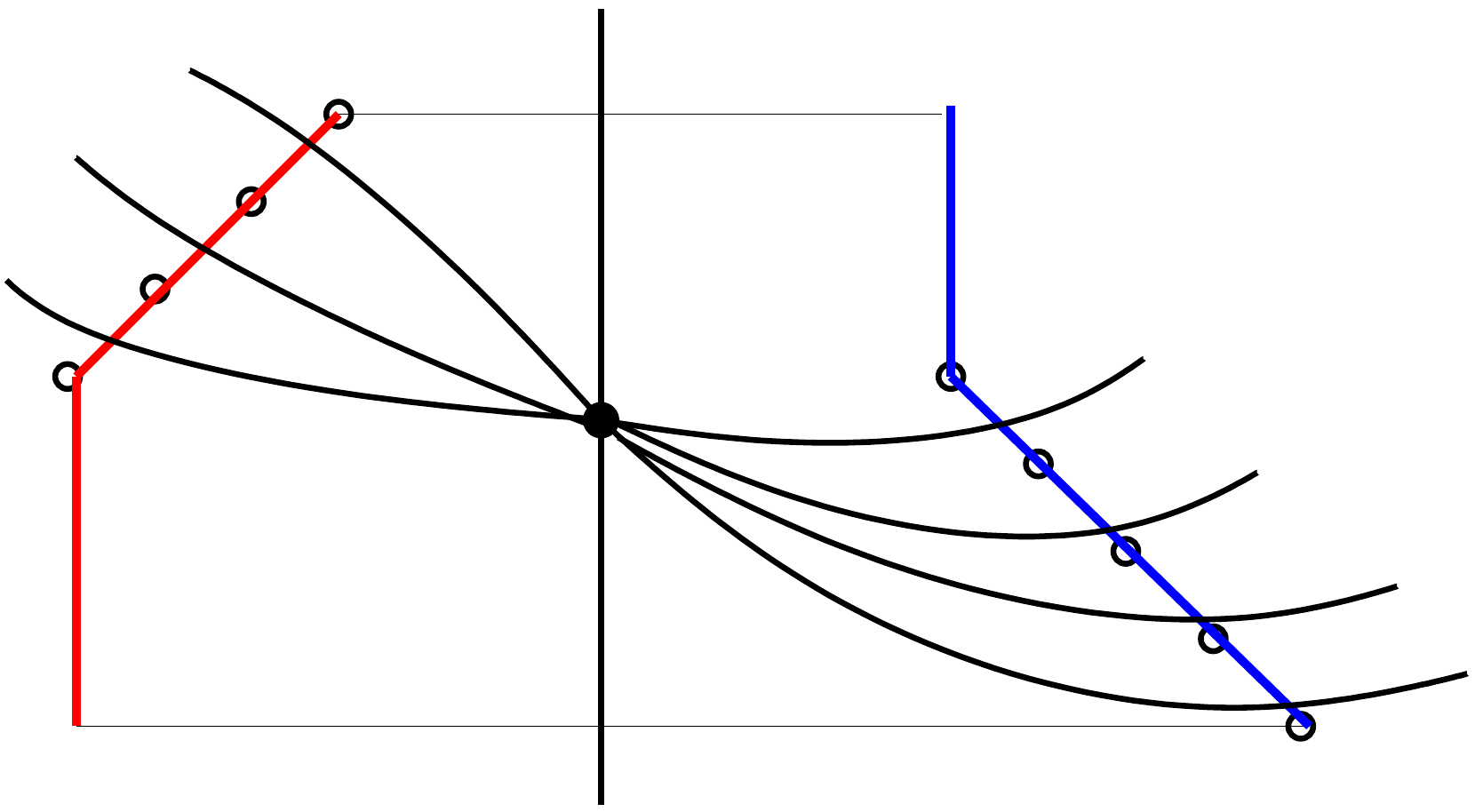}
  $$
There are four edges on the right, three on the left, and two horizontal ones (remember that vertical edges are contracted into points when we contract the paths). Observe that one can change the order in which the $|i|+|j|$ steps are made without changing the resulting map.

One has to be careful that there is a potential conflict between such rules and the contracting rules for oblique paths. We will nevertheless see some interesting examples.
\subsubsection{Other contractions}
In this paper I consider only contractions of quadrilaterals, according to one of the two non-crossing pairings of its sides. It would be possible to consider also contractions of larger faces, as constructed in section \ref{step}. Such contractions would be obtained from pairings of faces with an even number of sides. With non-crossing pairings one obtains planar maps, but it would be also possible to construct maps of higher genus by making more general pairings. However, we will not explore such possibilities in this paper.
\subsubsection{Pattern avoiding}
If we consider the allowed steps of the walk as forming an alphabet then the walk can be considered as  a word on this alphabet. It is possible to consider families of walks that are constrained to avoid certain patterns. Again we will encounter such examples. 
\subsubsection{Symmetries}There are some obvious symmetries in this construction. For example one can make a symmetry with respect to a vertical  axis in pictures like Figure \ref{fig:dual}. This corresponds to exchanging the horizontal and vertical directions of the walks, in other words to make a reflection with respect to the diagonal line $x=y$ in the plane.
Symmetry with respect to a horizontal  axis corresponds to considering walks with opposite steps, in reverse order.

\section{Maps with a spanning tree}\label{sec:span}
\subsection{}\label{sec:span1}
In this section we consider walks with steps in the set $\{(0,1),(-1,-1),(1,-1)\}$.
$${\tt    \setlength{\unitlength}{1.5pt}
\begin{picture}(30,20)
\thicklines

\put(0,10){\vector(1,-1){10}}
\put(0,10){\vector(-1,-1){10}}
\put(0,10){\vector(0,1){10}}
\end{picture}}
$$
For convenience we will give a name to these walks and refer to them as ``reversed $Y$-walks'', or $rY$-walks (note that, for  typographical reasons, the walks with opposite steps deserve the name of $Y$-walks).

\begin{figure}
$$
\includegraphics[width=14cm]{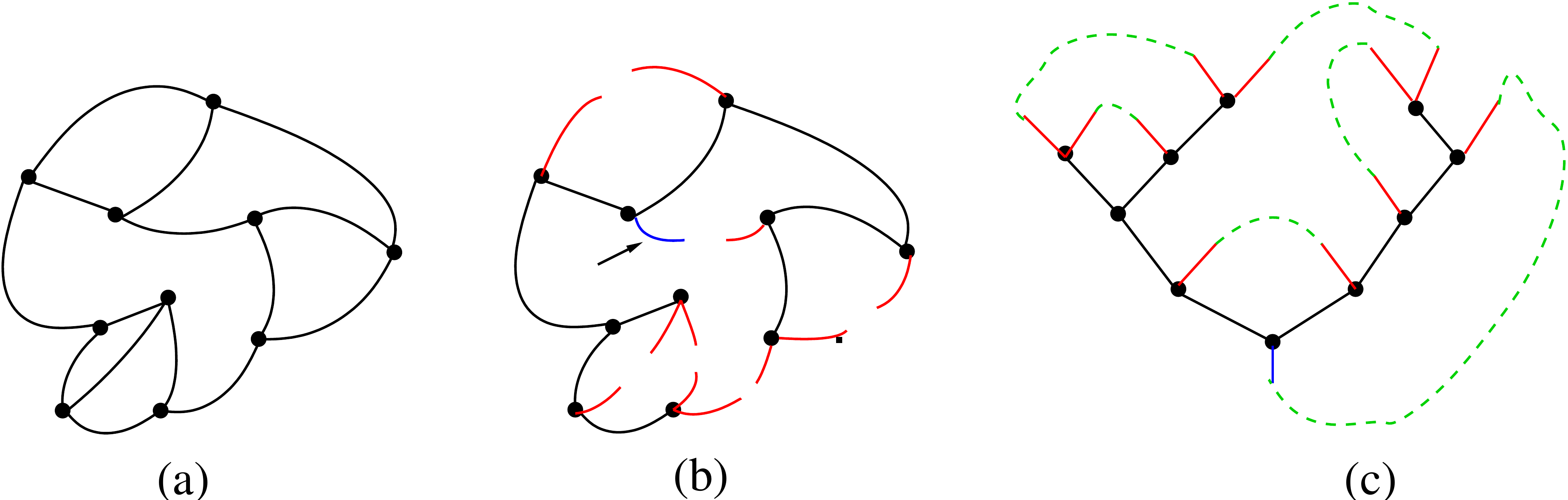}
 $$
\caption{A cubic map with a complete spanning tree and a marked leaf}
\label{fig:span}
\end{figure}

 It is easy to count the number of $rY$-walks in the quarterplane, starting and ending at $0$. Indeed the vertical coordinate of the walk gives a Dyck path, since the vertical coordinates of steps are either $1$ or $-1$. Choose a Dyck path of length $4n$  for the vertical coordinate and let $i_1,\ldots,i_{2n}$ be the indices of the down steps. The horizontal coordinate  moves only when a down step in the vertical direction is made, therefore the walk 
is specified by choosing another Dyck path of length $2n$, which describes the horizontal moves at times $i_1,\ldots,i_{2n}$.
It follows that the number of such walks is $C_{2n}C_{n}$ where $C_{l}=\frac{1}{l+1}{2l\choose l}$, the number of Dyck paths of length $2l$,
a  Catalan number. There is a simple model of maps which is counted by this number. 
Consider a cubic planar map (i.e. a map with vertices of degree three) with $2n$ vertices and thus $3n$ edges.
Take $n+1$ of the edges and cut them in two, in order to obtain a tree, with the half-edges being the leaves of the tree. We call this tree a ``complete'' spanning tree. If we choose one of these leaves to root the tree we get   a planar binary tree with $2n$ internal vertices and $2n+1$ leaves, plus the root leaf.
This construction is reminiscent of so-called ``blossoming bijections'' (see \cite{Sch97}, or \cite{AP} for a recent reference), except 
that we do not orient the leaves of the tree.
 The map from which we started can be recovered by pairing these $2n+2$ leaves, to reconstruct the edges which have been cut. In Figure \ref{fig:span}(a) I show  a cubic map and a complete spanning tree in \ref{fig:span}(b). The root half-edge is in blue, with an arrow pointing at it  and the other half-edges are in red.  In Figure \ref{fig:span}(c)  I show the planar binary tree in standard representation, with the matching of the leaves giving back the original map.

Such objects, planar cubic maps on $2n$ vertices, with a complete spanning tree with a marked leaf, are thus in bijection with pairs $(T,P)$ where $T$ is a binary tree with $2n$ internal vertices and a leaf added to the root, and $P$ is a pairing of the $2n+2$ leaves. The leaves can be ordered cyclically from the root by going clockwise around the tree, thus the pairings are in bijection with non-crossing pairings 
of $[1,2n+2]$.
There are $C_{2n}$ planar binary  trees. Since there are $C_{n+1}$ non-crossing pairings we get $C_{2n}C_{n+1}$ pairs $(T,P)$. This is not quite $C_{2n}C_{n}$, but if we add the constraint that the root leaf has to pair with its immediate successor then we get the right number $C_{2n}C_{n}$. 
 We call such an object a ``special'' planar cubic map on $2n$ vertices, with a complete spanning tree with a marked leaf.
For example the cubic map depicted in Figure \ref{fig:span} is special since the root leaf is paired with its successor.
  We will see that our construction provides a geometric realization of this bijection between walks and such maps. We will also extend it to  a walk model for ``non-special'' cubic maps with a complete spanning tree rooted at some half-edge.

\subsection{The construction}
The $rY$-walks have two types of oblique steps, so we must give the rules for contracting the associated quadrilaterals. These rules are simple, we show them below.
\begin{contraction}\label{contract:1}
$${\tt    \setlength{\unitlength}{1pt}
\begin{picture}(300,65)
\thicklines    
\
\put(-13,-10){$w$}\put(48,-10){$z$}\put(-1,16){$u$}\put(68,16){$v$}
\put(-5,0){\line(1,0){60}}
\put(55,0){\color{blue}\line(1,1){10}}
\put(5,10){\line(1,0){60}}

\put(-5,0){\color{red}\line(1,1){10}}           
 \put(80,3){$\to$}
\put(98,3){$w$}
\put(110,5){\color{red}\line(1,0){40}}\put(150,5){\color{blue}\line(1,0){40}}\put(148,8){$u$}\put(148,-4){$z$}\put(195,3){$v$}
\thicklines    
\
\put(-8,32){$w$}\put(63,32){$z$}\put(-16,56){$u$}\put(73,56){$v$}
\put(0,40){\line(1,0){60}}
\put(60,40){\color{blue}\line(1,1){10}}
\put(-10,50){\line(1,0){80}}

\put(0,40){\color{red}\line(-1,1){10}}           
 \put(80,43){$\to$}
\put(100,43){$w$}
\put(110,45){\color{red}\line(1,0){40}}\put(150,45){\color{blue}\line(1,0){40}}\put(148,48){$u$}\put(148,36){$z$}\put(193,43){$v$}
\end{picture}}
$$
\end{contraction}

\medskip

i.e. in both cases  identify the NW and SE corners.
\subsection{Some bijections}
\begin{figure}
\caption{$rY$-walks and their associated maps.\newline (a) Construction of the spanning tree using the vertical coordinates.\newline (b) Matching the leaves using the horizontal coordinates.}
\label{rY-walks}
\includegraphics[width=10cm]{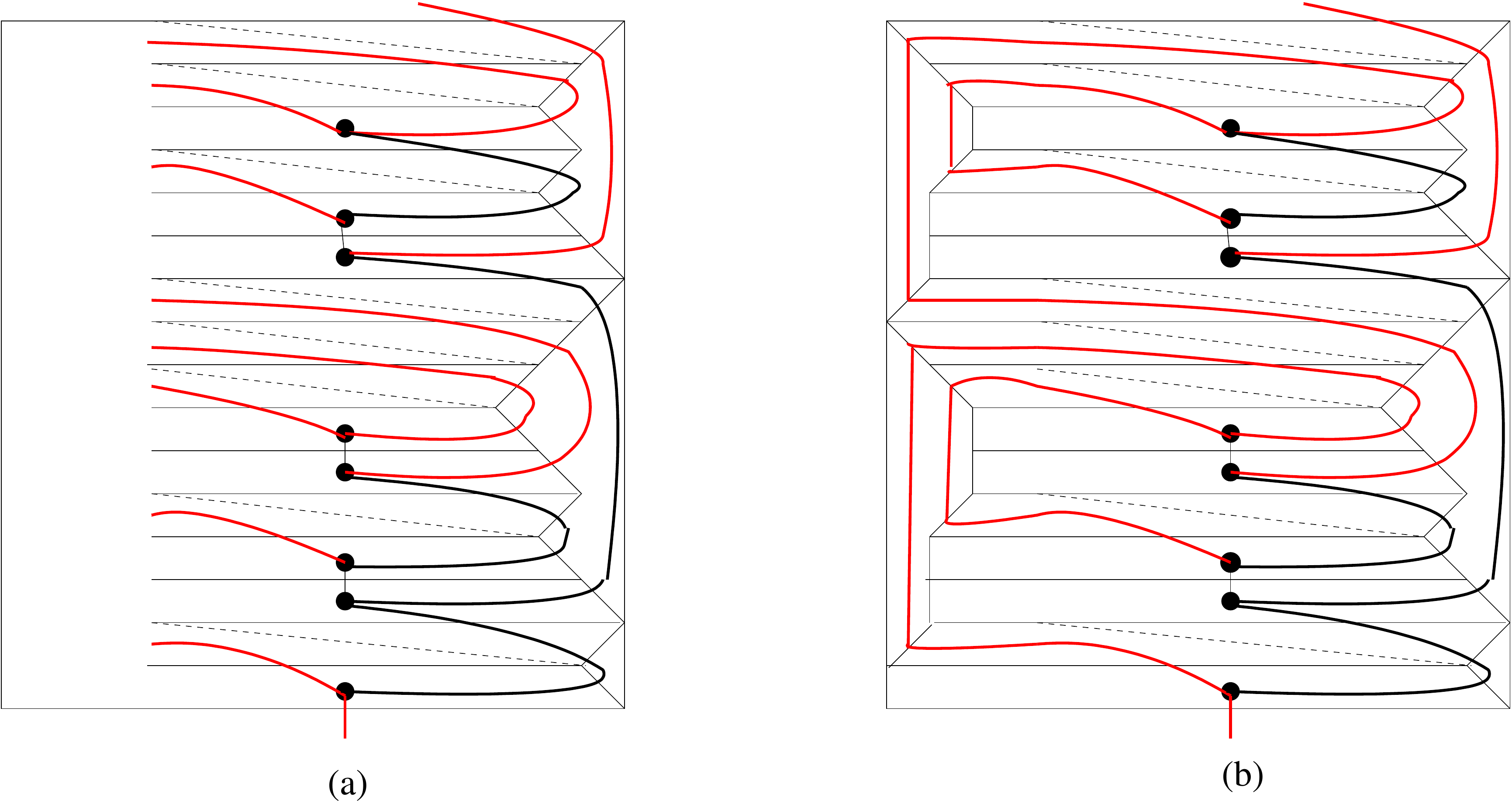}
\end{figure}

There is a simple and well known bijection, which we shall denote by $\mathcal T$,  between rooted planar binary trees with $n$ vertices and Dyck paths with $2n$ steps, which is best described recursively. Write a Dyck path as a sequence of up $(u)$ and down $(d)$ steps. Then $\mathcal T(ud)$ is the unique rooted planar binary tree with one vertex. For   any Dyck path
$\mathcal{ D}$,  of length $2n$,  write it uniquely as $\mathcal{ D}=u\mathcal{ D}_1d\mathcal{ D}_2$, where $\mathcal{ D}_1$ and $\mathcal{ D}_2$ are Dyck paths of lengths $2k-2$ and  $2n-2k$. Here $2k$ is the time of first return to $0$ for $\mathcal{ D}$. Then  $\mathcal T(\mathcal{ D})$ is the binary tree whose left branch is $\mathcal T(\mathcal{ D}_1)$ and whose right branch is $\mathcal T(\mathcal{ D}_2)$.
This tree can be obtained easily from our construction as follows: write the Dyck path on the right but leave the left extremities of the quadrilaterals unfinished, so that we do not make the identifications on the left of the picture. This is shown in Figure \ref{rY-walks}(a).
The cut dual map is a tree and it is obvious that its construction satisfies the same recursion as the bijection $\mathcal T$.
Observe that the leaves of the tree (in red in the picture) appear ordered on the left of the picture (except the root leaf and its immediate successor).
There is also a simple and well known bijection between Dyck paths and noncrossing pairings which consists in pairing up and down steps. For example, the  non-crossing pairing $(1,8),(2,3),(4,7),(5,6)$ corresponds to the Dyck path $uuduuddd$;
$$
\includegraphics[width=4cm]{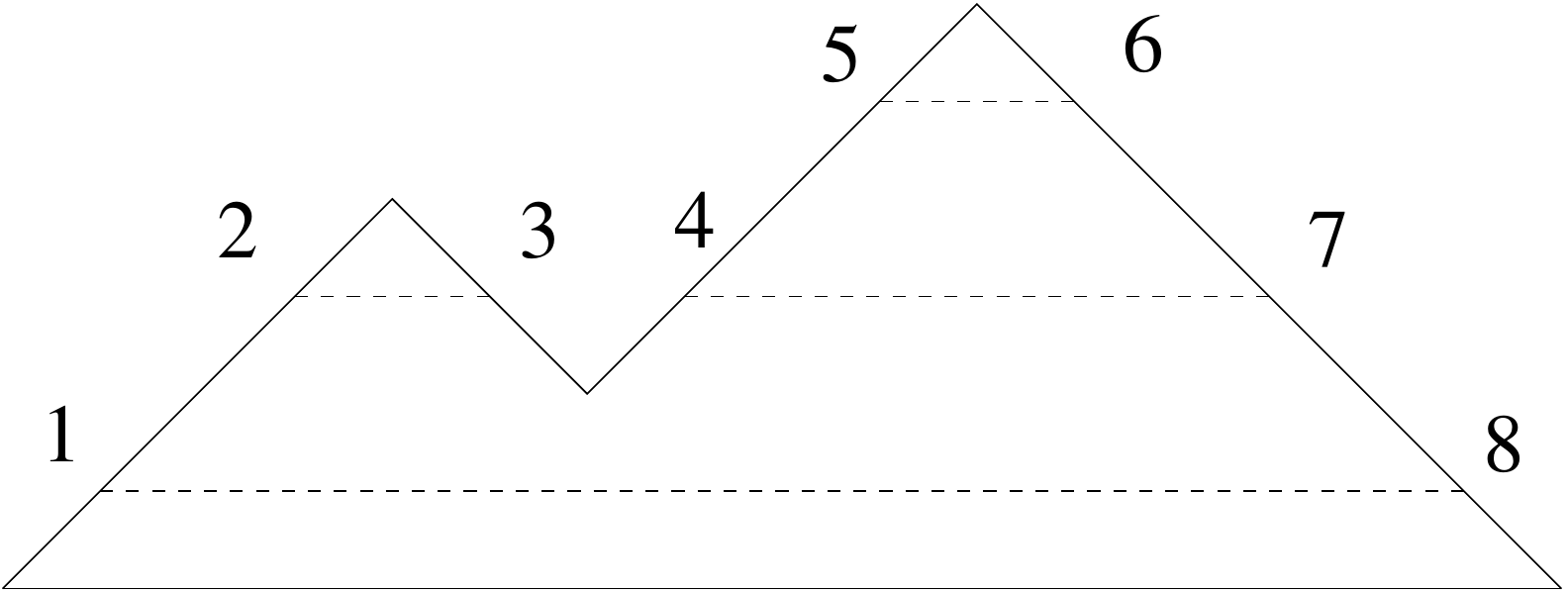}
 $$

\begin{figure}
\caption{Non special maps.\newline (a) Map with four unpaired leaves.\newline (b) Matching the remaining leaves .}
\label{rY-walks-sp}
\includegraphics[width=10cm]{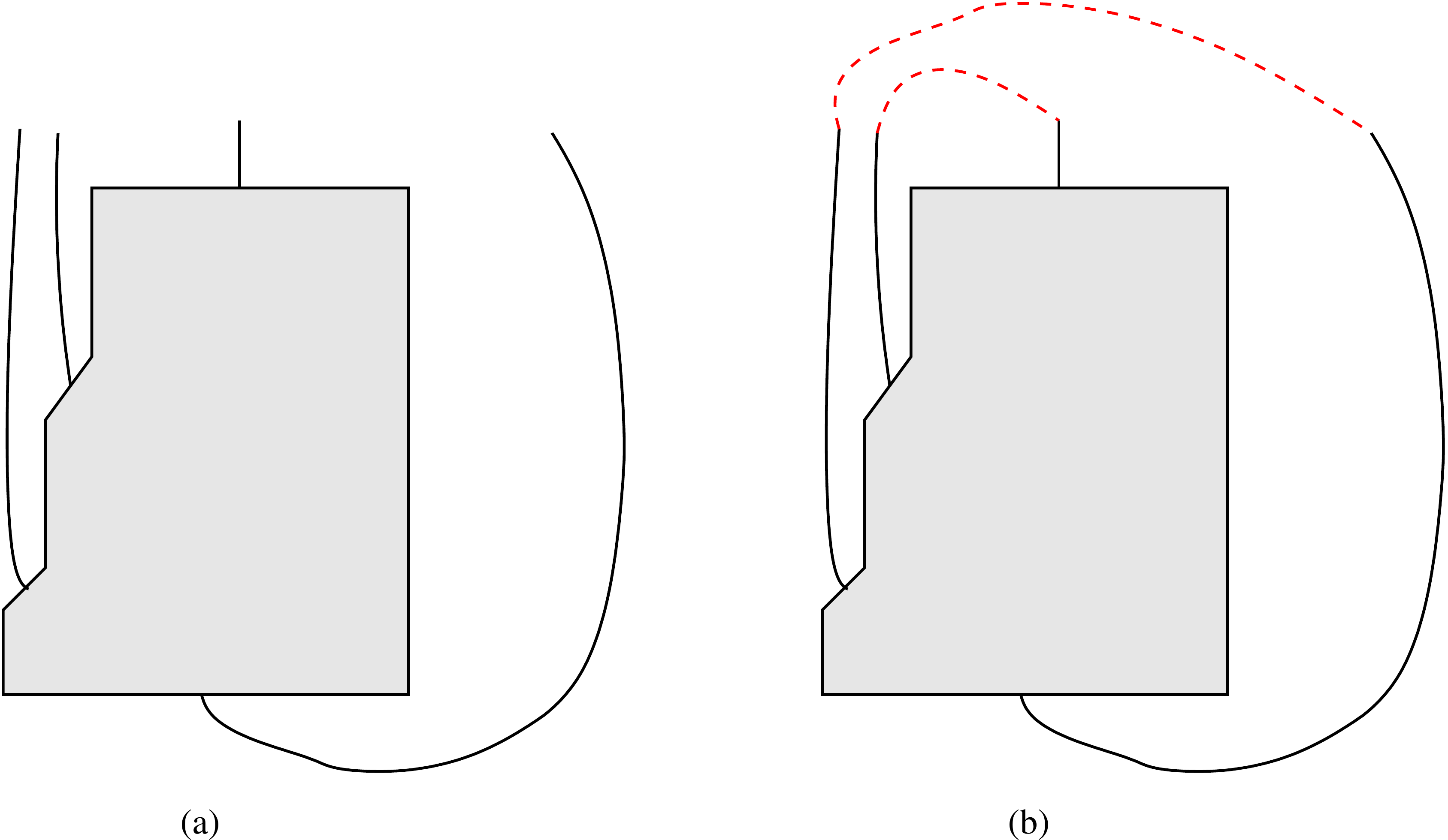}
\end{figure}

When we complete the picture on the left, as in Figure \ref{rY-walks}(b)
the Dyck path corresponding to the horizontal coordinate determines the non-crossing pairing of the leaves.
It  remains to pair the root leaf with the leaf going through the upper side of the rectangle to obtain a special cubic map. Conversely, given a special cubic map with a complete, rooted  spanning tree, the spanning tree gives us the vertical coordinate of an $rY$-walk, while the pairing of the leaves gives us the horizontal coordinate. Finally one has the following.
\begin{theorem}
The construction with Contraction Rules \ref{contract:1} gives a bijection between  $rY$-walks in the quarter plane, starting and ending at $0$, with $4n$ steps,
and special cubic maps with a complete spanning tree, rooted at some half-edge,  with $2n$ vertices. 
\end{theorem}
We have seen that, in this construction,  the root leaf is paired with its immediate successor. In order to obtain all possible pairings, let us instead consider walks starting at $0$ but ending at the point $(2,0)$. Making the contractions of paths as in section \ref{sec2:2}
we end up with Figure \ref{rY-walks-sp}(a), where there are four half-edges unmatched in the dual map: the half-edges corresponding to the upper and lower sides of the rectangle and two half-edges coming from the two up steps on the left which are not matched with a down step.
We then pair these half-edges so that the upper and lower horizontal sides are not matched, as in Figure \ref{rY-walks-sp}(b)
This gives a cubic map with a complete spanning tree, in which the root is not paired with its immediate successor. 

Putting this together with the preceding construction, we get:

\begin{theorem}
The construction with Contraction Rules \ref{contract:1} gives a bijection between $rY$-walks on the quarter plane  starting  at $0$ and ending at $0$ or $(2,0)$, with $4n$ steps,
and  cubic maps with a rooted complete spanning tree, with $2n$ vertices. 
\end{theorem}

\subsection {Variants}
\begin{figure}
\caption{Quartic maps.}
\label{Fuss-walks}
\includegraphics[width=10cm]{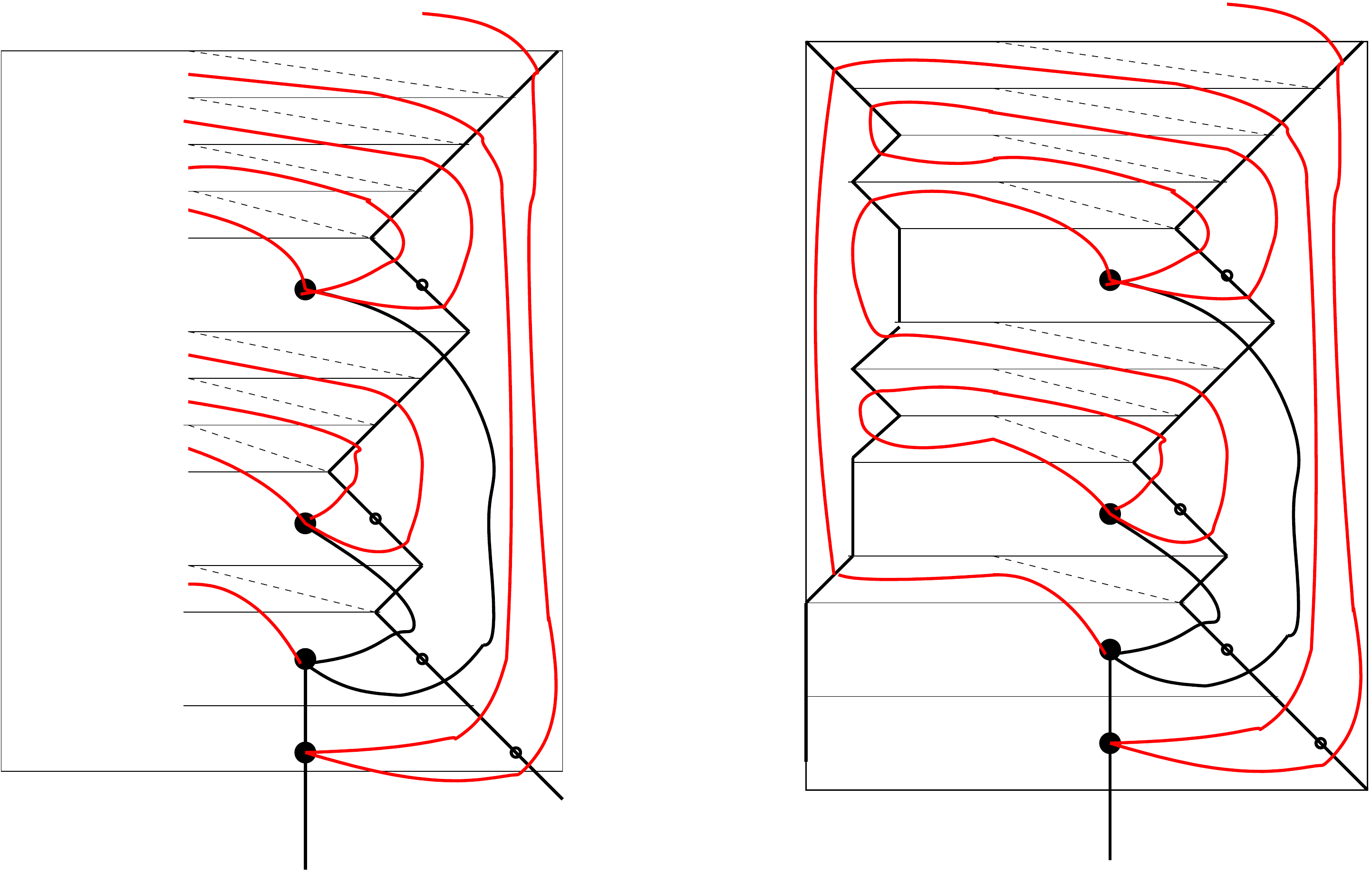}
\end{figure}

Instead of the step $(0,1)$ we could consider  $(0,2)$. Again it is easy to enumerate all walks with step set $\{(-1,-1),(1,-1),(0,2)\}$ which start and end at $0$: in the vertical direction we have a walk with steps $2$ and $-1$ and an application of the cycle lemma give the number of such walks with $3n$ steps, which is $\frac{1}{2n+1}{3n\choose n}=C^{(2)}_n$, a Fuss-Catalan number (see e.g. \cite{Bernardi}). On the set of $2n$ negative steps we have to choose a Dyck path of length $2$ which gives a Catalan number $C_n$. The final count is 
$C^{(2)}_nC_n$. Again our construction  gives a bijection between walks with step set $\{-1,2\}$ and ternary trees, as in Figure 
\ref{Fuss-walks}(a), which we can complete as in  Figure 
\ref{Fuss-walks}(b).
\begin{theorem}
The construction with Contraction Rules \ref{contract:1} gives a bijection between walks on the quarter plane with step set $$\{(0,2),(-1,-1),(1,-1)\}$$ starting and ending at $0$, with $3n$ steps,
and special quartic maps with a complete spanning tree, with $3n$ vertices. 
When extending to 
 walks ending at $(2,0)$, we get
a bijection with quartic maps with a complete rooted spanning tree. 
\end{theorem}

 Obviously this can be further extended to a bijection between  planar maps equipped with a spanning tree and a marked leaf,  with walks whose steps belong to the  infinite  set
$$\{(-1,-1),(1,-1),(0,k);k\geq 0\}$$ starting at $0$ and ending at $0$ or $(2,0)$. Observe that the vertical component of such a walk is a
\L{}ukasiewicz path. Applying the same algorithm as above recovers a well known  bijection between \L{}ukasiewicz  paths and trees   (see e.g. \cite{Stanley}), which extends the bijection between Dyck paths and binary trees. 
\begin{theorem}
The construction with Contraction Rules \ref{contract:1} gives a bijection between walks on the quarter plane with step set 
$$\{(-1,-1),(1,-1),(0,k),k\geq 0\}$$ starting at $0$ and ending at $0$ or $(2,0)$, and planar maps with vertices of degree $\geq 2$, 
with a complete spanning tree. In this bijection, steps of type $(0,k)$  of the walk correspond to vertices of degree $k+2$ of the map.
\end{theorem}
 Mullin's construction (as explained in Schaeffer \cite{Schaeffer}) also provides a bijection between maps equipped with a spanning tree and walks in the quarter plane but this is a quite different bijection, for example the set of steps allowed in Mullin's bijection is the set of straight steps.
\section{Recovering other bijections}\label{sec:4}
\subsection{Triangulations with a Hamiltonian cycle on the faces}
We consider a planar triangulation equipped with a Hamiltonian cycle of its dual graph. This amounts to choosing an  ordering $f_1,f_2,\ldots,f_n$ of the faces of the triangulation  such that, for all $i$, the faces $f_i$ and $f_{i+1}$ are adjacent (where $i+1$ is taken modulo $n$) and for each $i$ one chooses an  edge $e_i$ in the common boundary of $f_i$ and $f_{i+1}$. Here is an example where the Hamiltonian cycle is in red, the faces are numbered from $1$ to $8$ and the edges $e_i$ are the edges crossed by the red path.
$$
\includegraphics[width=6cm]{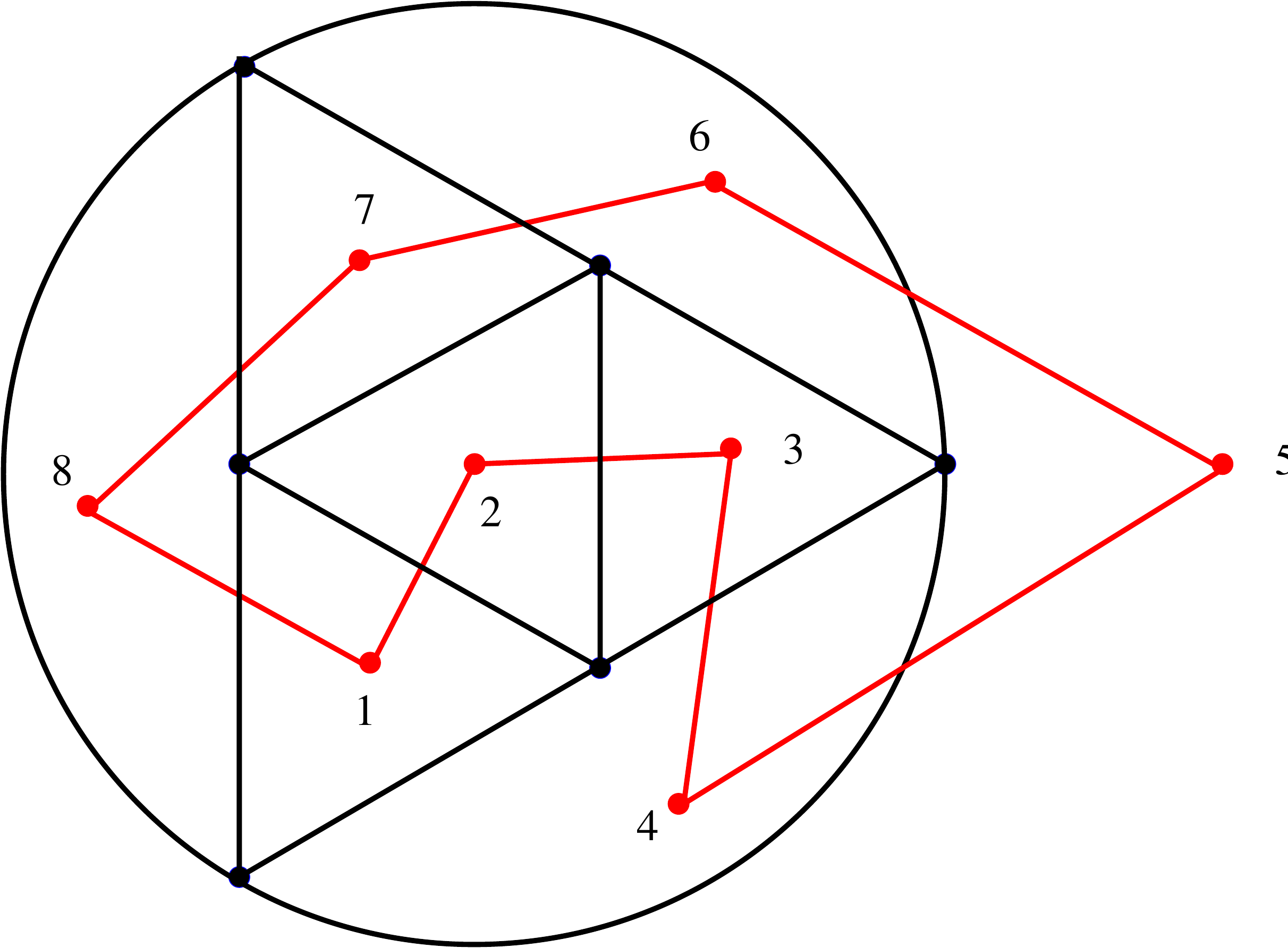}
  $$

The graph formed from the vertices of the triangulation and the edges which do not belong to the set $e_1,e_2,\ldots,e_n$ is made of two disjoint trees and the Hamiltonian cycle goes around each of them, one of them being on its right and the other one on its left. Here is the picture:

$$
\includegraphics[width=6cm]{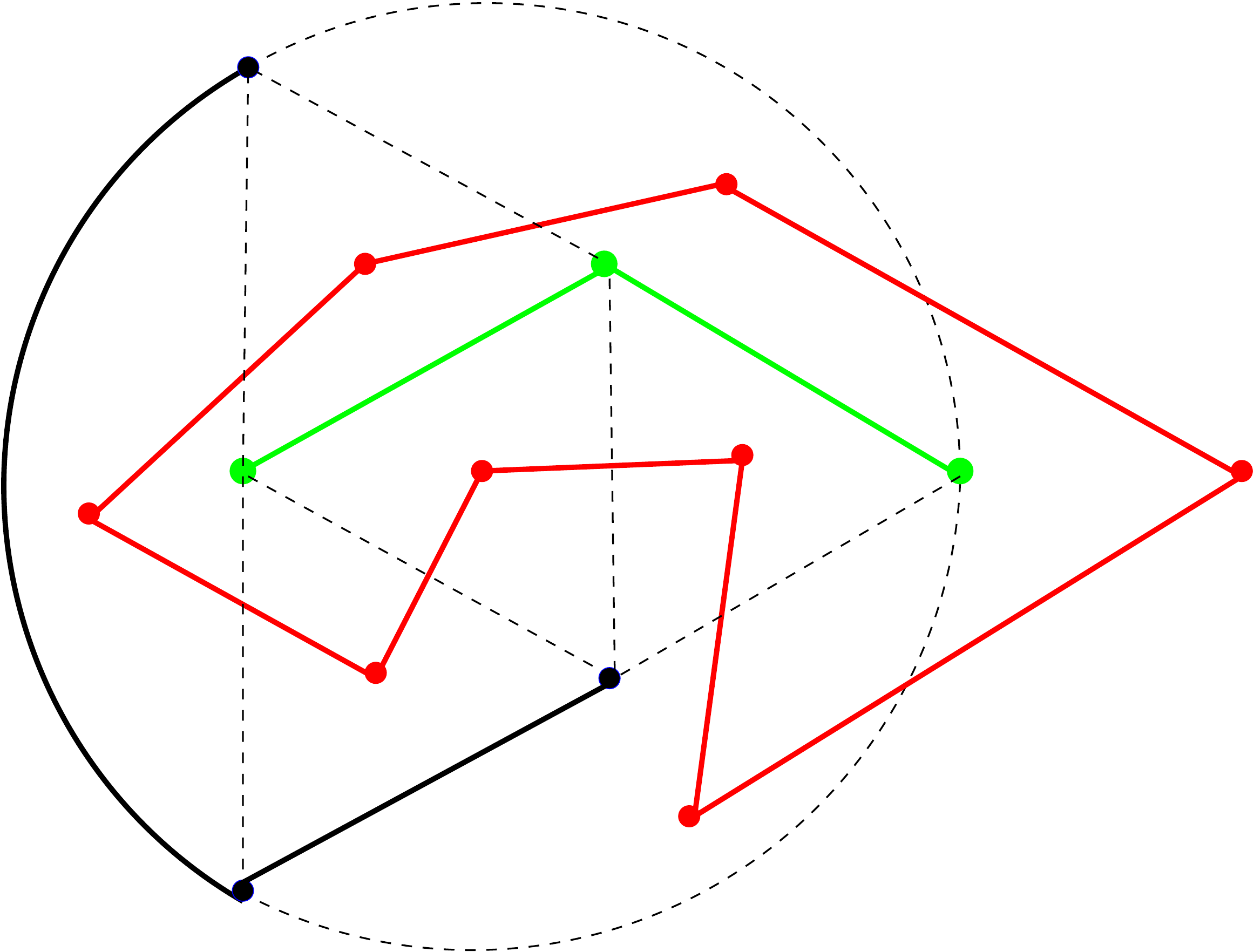}
  $$

We go around the cycle  and record the successive triangles. Each such triangle has exactly one side in one of the trees.
If this side is in the left tree  we record a $(1,0)$ step,   if the path  goes up  in the tree and a $(-1,0)$ step if it goes down. If the side is in the right tree we record  similarly a $(0,\pm 1)$ step. The final picture is 
$$
\includegraphics[width=4cm]{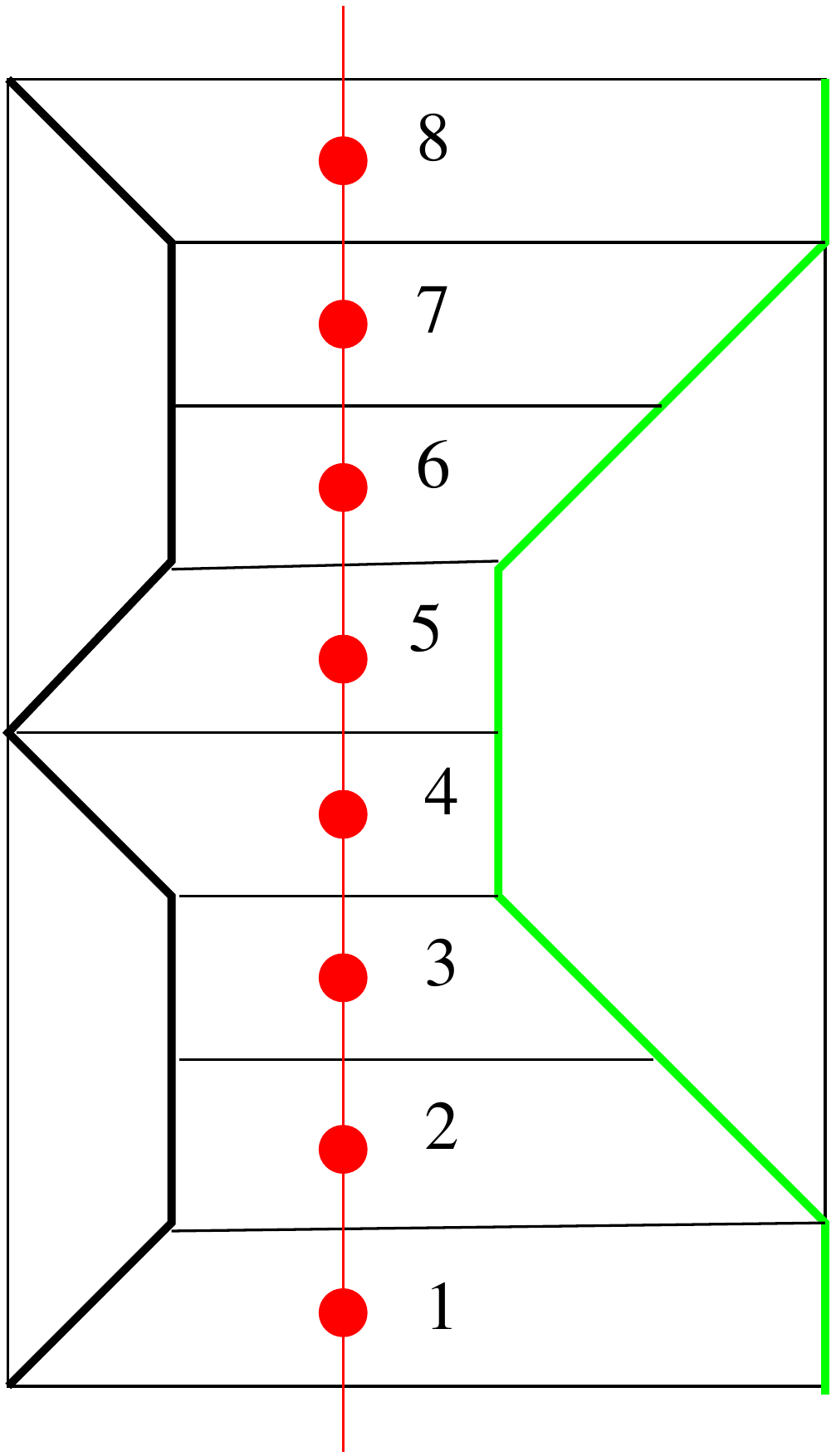}
  $$
where we draw the Hamiltonian cycle in red, but not the whole dual map.
It is immediate to see that we obtain in this way a walk in the quarter plane, with straight steps, starting and ending at $0$. This construction  yields a bijection between walks with straight steps in the quarter plane, starting and ending at $0$, and
planar triangulations  equipped with a Hamiltonian cycle of the dual graph.
 This is well known and occurs in Mullin's bijection between maps with a spanning tree and walks on the quarter plane, see e.g. \cite{Schaeffer}, section 1.2. See also \cite{GHS} for a recent overview of applications of this kind of bijections to random maps.
In the remaining sections we will obtain examples with more flexibility by looking at classes of walks having oblique steps.
\subsection{ Kreweras walks and  Bernardi's  bijection }
We now consider a bijection of Bernardi \cite{Bernardi} (also used in \cite{Bernardi-Holden-Sun}). 
A Kreweras walk has steps in the set  $\{a,b,c\}$ where $a=(1,0),$ $b=(0,1)$ and $c=(-1,-1)$. 

$${\tt    \setlength{\unitlength}{1.5pt}
\begin{picture}(30,20)
\thicklines

\put(0,10){\line(1,0){10}}\put(12,9){\text{$a$}}

\put(0,10){\line(-1,-1){10}}\put(-16,-2){\text{$c$}}

\put(0,10){\line(0,1){10}}\put(-1,22){\text{$b$}}
\end{picture}}
$$

It can be shown that Kreweras walks with $3n$ steps, which remain in the quarter plane, starting and ending at $0$, are enumerated by the formula $\frac{2^n}{2n+1}{3n\choose n}$. This  formula is explained in a bijective way in  Bernardi \cite{Bernardi}. We will see  how to recover his bijection using our construction. 

We thus consider a Kreweras walk, which remains in the quarter plane, starting and ending at $0$,
When we do the construction of section \ref{sec:3.1} the steps of types  $a$ or $b$ give triangles when we contract the Motzkin paths. The steps of type 
 $c$ give rise to quadrilaterals as below:
$${\tt    \setlength{\unitlength}{1pt}
\begin{picture}(300,20)
\thinlines    
\
\put(-7,-7){$w$}\put(63,-6){$z$}\put(-16,16){$u$}\put(73,16){$v$}
\put(0,0){\line(1,0){60}}
\put(60,0){\color{blue}\line(1,1){10}}
\put(-10,10){\line(1,0){80}}

\put(0,0){\color{red}\line(-1,1){10}}           
 
\end{picture}}
$$

\medskip

\begin{contraction}
Consider the  sides   $uw$ and $vz$ of the quadrilateral corresponding to a step of type $c$. In the contraction of the two Motzkin paths, each of these sides is matched with another  segment
below it, moreover this segment belongs to a step of type $a$ (for the $uw$ segment) and to a step of type $b$ (for the $vz$ segment).
In the enumeration $x_1,x_2,\ldots,x_n$ of the steps of the walk let $i$ and $j$, respectively, be  the indices of these steps, thus
$x_i=a$ and $x_j=b$ . If $i<j$ then we identify  $u$ and $z$ in the contraction of the quadrilateral $uvwz$. If  $i>j$ then we identify $v$ and $w$.

\end{contraction}

Here is the case  of the path with sequence of steps $aabbccbac$, with the dual cubic map in green: 
$$
\includegraphics[width=4cm]{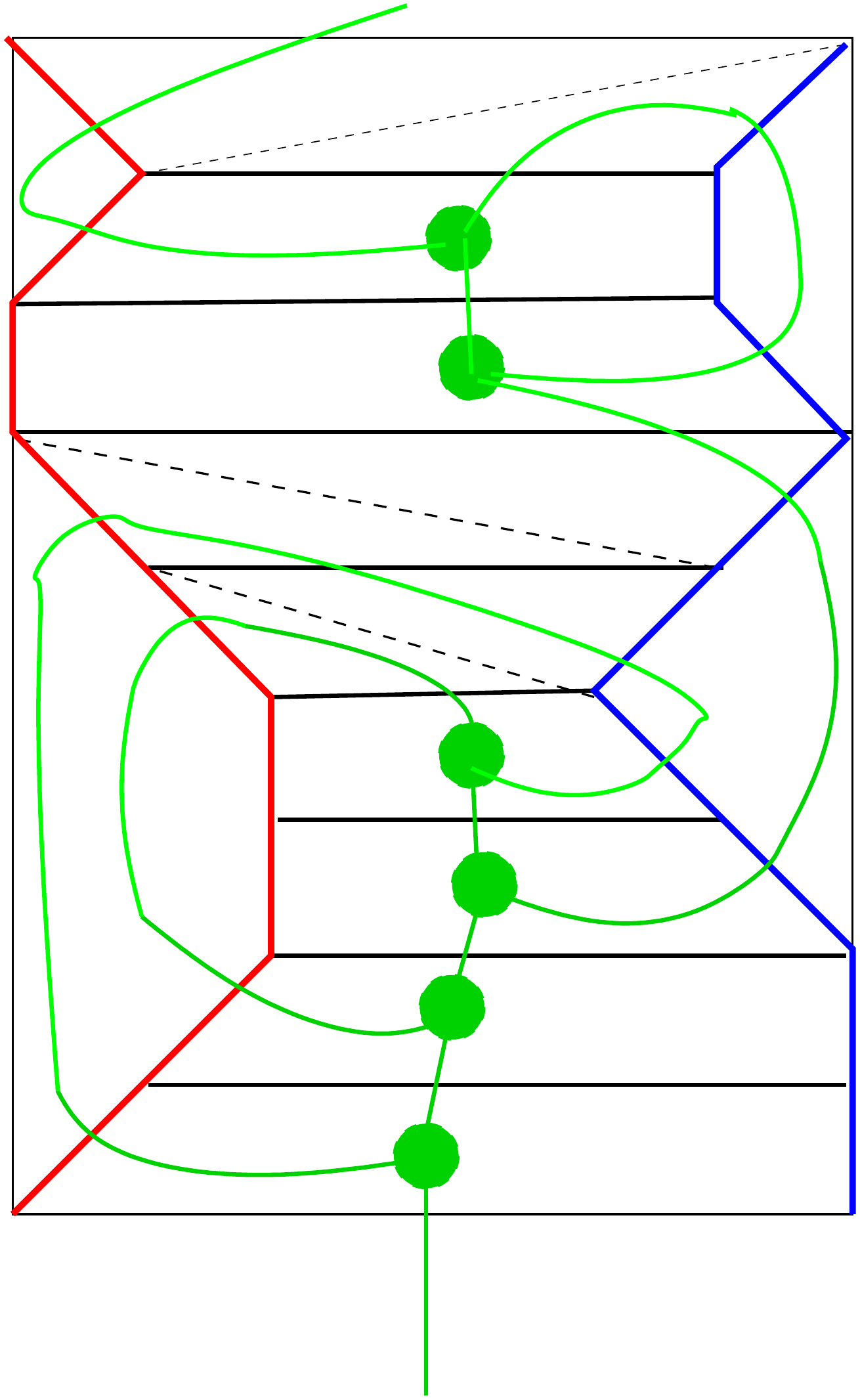}
 $$
 Observe that the triangulation (and the dual cubic map) are loopless. Indeed a loop in the construction could be obtained only if, just after an $a$ or $b$ step, the $c$ step is contracted so that it has two sides in common with the preceding step:

$$
\includegraphics[width=4cm]{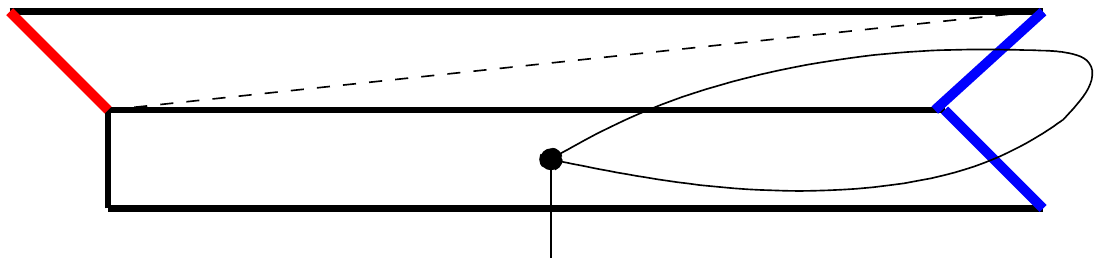}
 $$
but the rules of contraction prevent this.

Each vertex of the dual map has three half-edges, which we orient so that the half-edge corresponding to the base of the triangle is incoming and the ones corresponding to the other sides of the triangle are outgoing:

$${\tt    \setlength{\unitlength}{1pt}
\begin{picture}(300,40)
\thicklines    
\put(25,18){\color{green}\circle*{4}}
\put(25,18){\color{green}\vector(-3,2){18}}
\put(25,18){\color{green}\vector(2,3){12}}
\put(25,0){\color{green}\vector(0,1){16}}

\thinlines  
\put(0,10){\line(1,0){60}}\put(0,10){\line(1,1){20}}\put(20,30){\line(2,-1){40}}

\end{picture}}
$$

If we keep only the edges of the dual map having half-edges with matching orientations we obtain a spanning tree like this:
$$
\includegraphics[width=4cm]{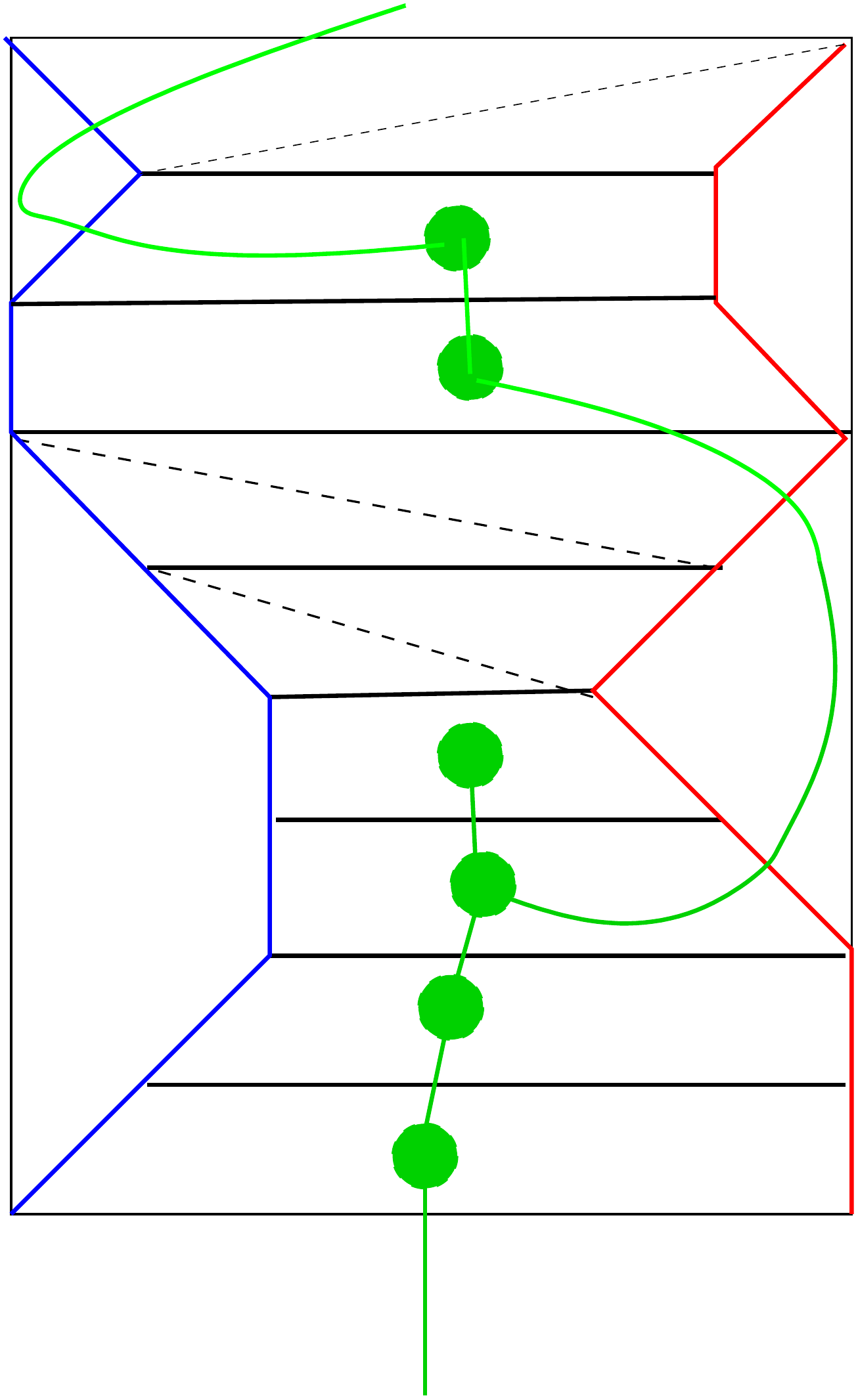}
 $$
We will now see that the data of the loopless triangulation and the spanning tree of the dual are exactly the ones obtained from Bernardi's bijection. The tree is a special kind of tree, called depth tree in \cite{Bernardi}, but we will not need to go into this here.
Let us sketch Bernardi's construction and check that it is the same as this one. More details can be found in \cite{Bernardi}. 
The construction is done step by step, by growing the dual cubic map and a spanning tree, using three mappings denoted $\varphi_a,\varphi_b,\varphi_c$. Observe that in \cite{Bernardi} the Kreweras walk has in fact opposite steps as the ones used in this paper, but since Bernardi scans the steps of the walk in reverse order,  the result is the same.
The algorithm starts with a vertical arrow pointing outside the root:

$$
\includegraphics[height=1cm]{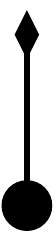}
 $$

 At each step a growing map (which we depict by a circle) is constructed, with  half-edges pointing outside, one of them being  endowed with an arrow pointing upwards:

$$
\includegraphics[height=1.7cm]{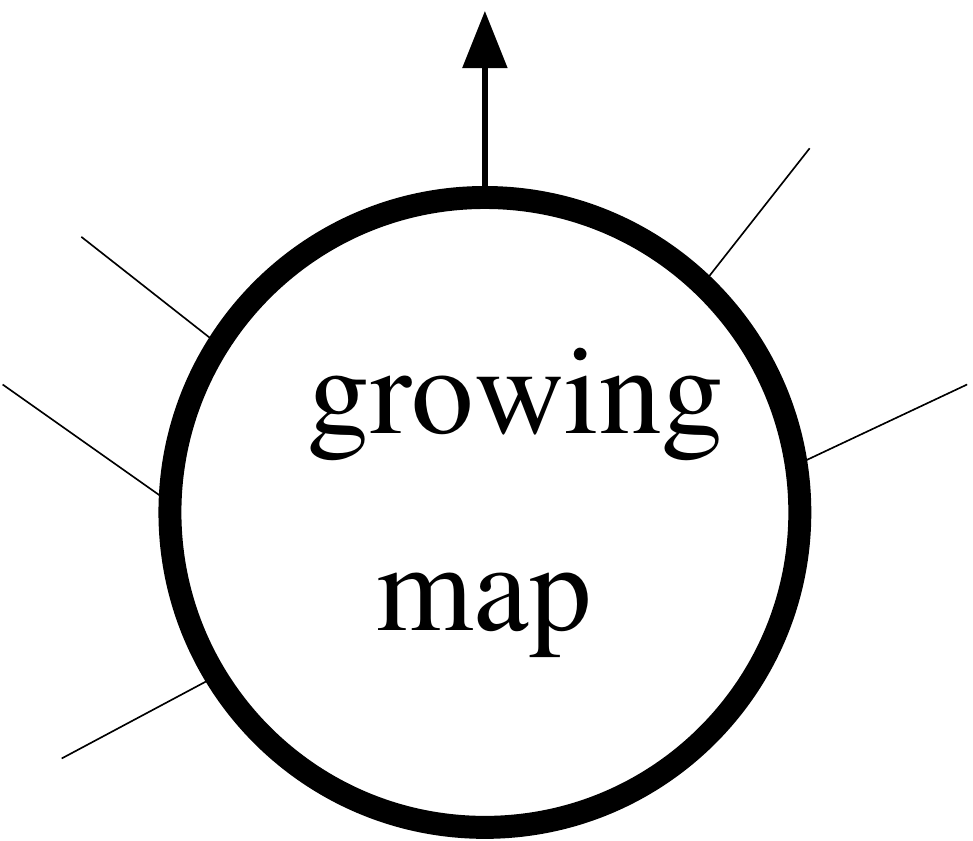}
 $$
The mapping $\varphi_a$  consists in transforming the half-edge containing the arrow into an edge and adding on top of this edge a pair of half-edges, with the right half-edge carrying the arrow like this:
$$
\includegraphics[height=2cm]{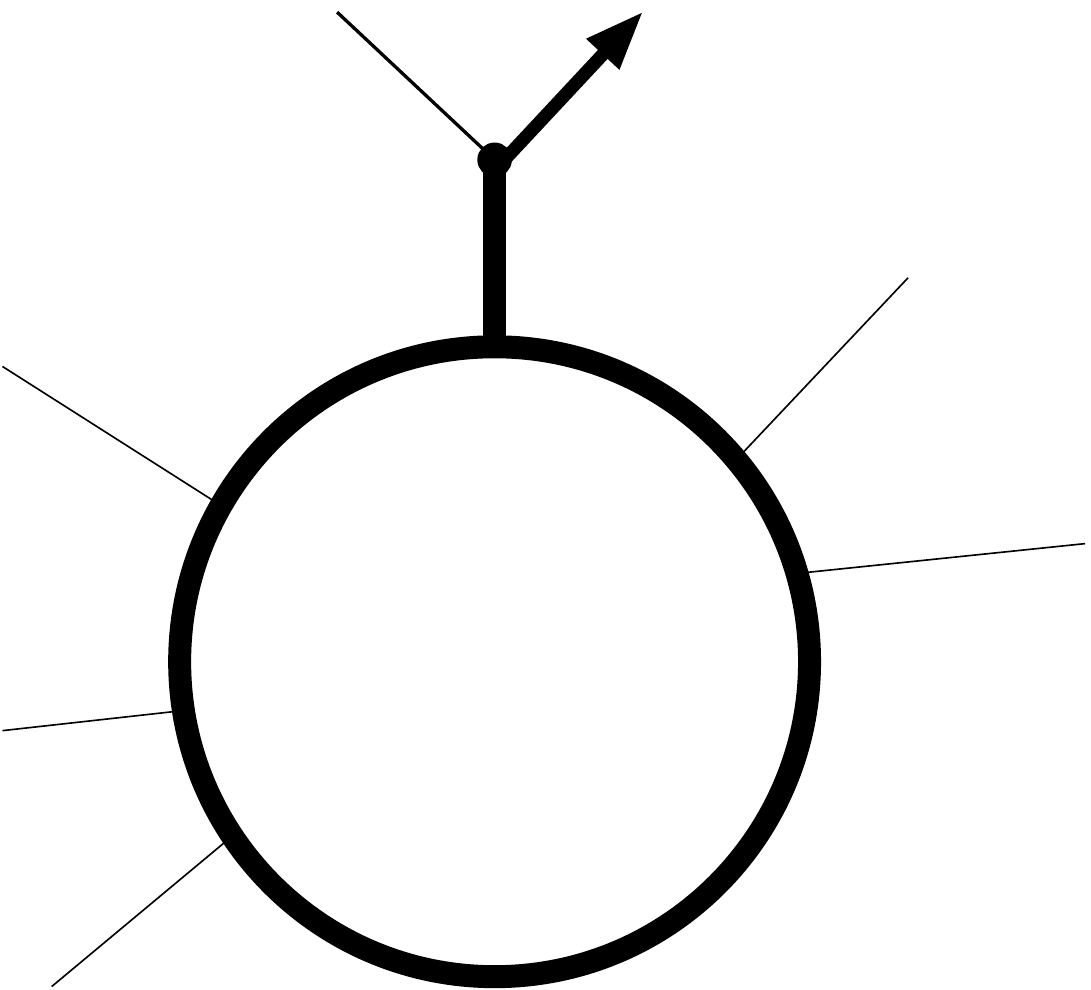}
 $$

In case of a step $b$ the mapping $\varphi_b$ is similar but the  arrow is on the left. 

$$
\includegraphics[height=2cm]{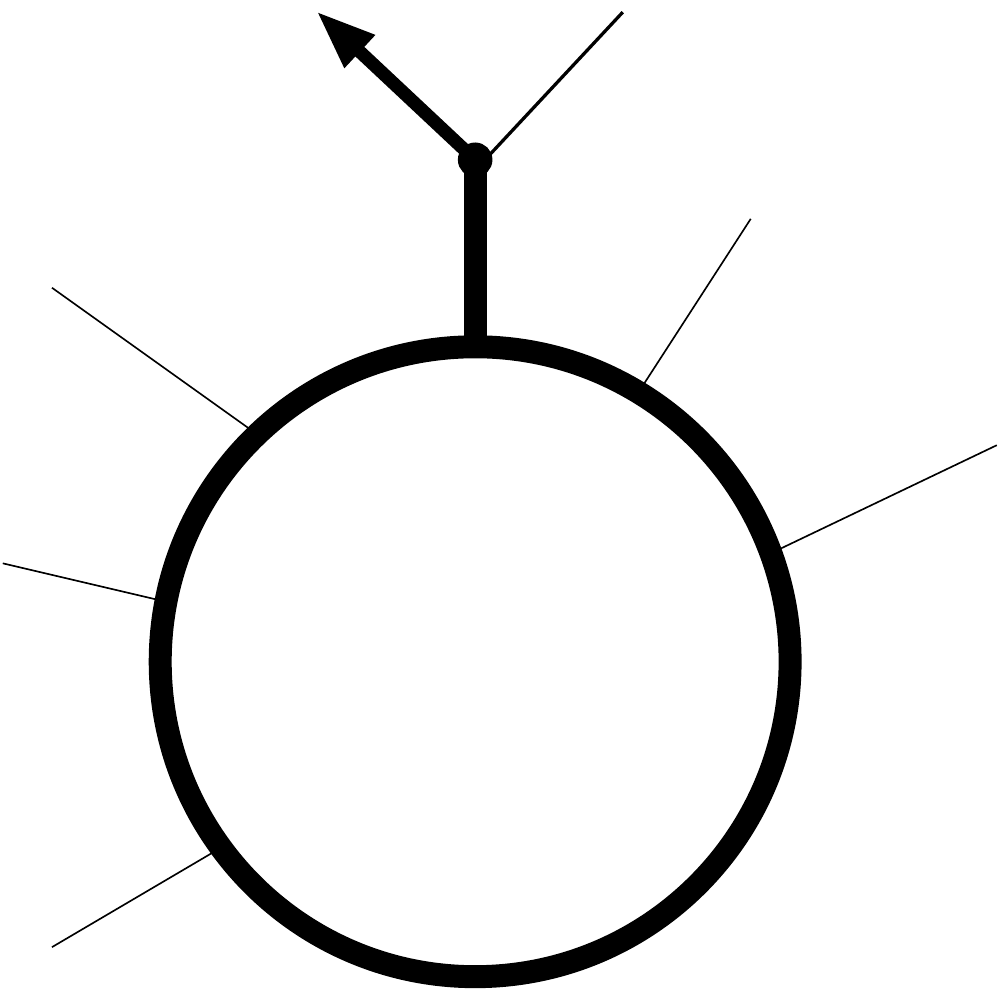}
 $$
Among all half-edges pointing out of the growing map there exists  some  on the right of the arrow and some on the left, in particular if there are at least one half-edge on each side we can single out the ones which are closest to the arrow. 

Bernardi shows that in his construction one of them is an ancestor of the other in the growing tree. Let us call $s$ this ancestor and $t$ the other half-edge.
$$
\includegraphics[height=2.6cm]{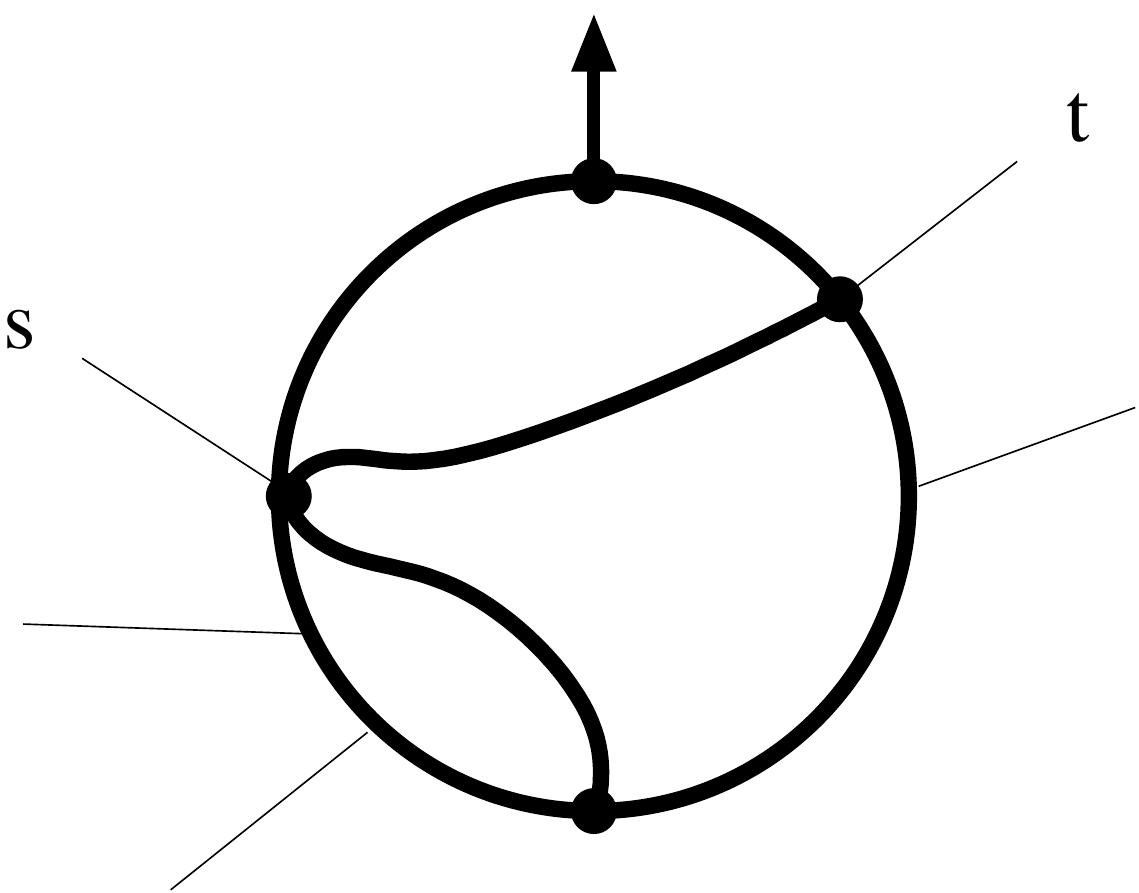}
 $$ 
In case of a step $c$ the mapping $\varphi_c$ consists  in joining the half-edge containing the arrow to $s$ in order to make an edge of the growing map and taking $t$ to carry  the new arrow. 
$$
\includegraphics[height=2.5cm]{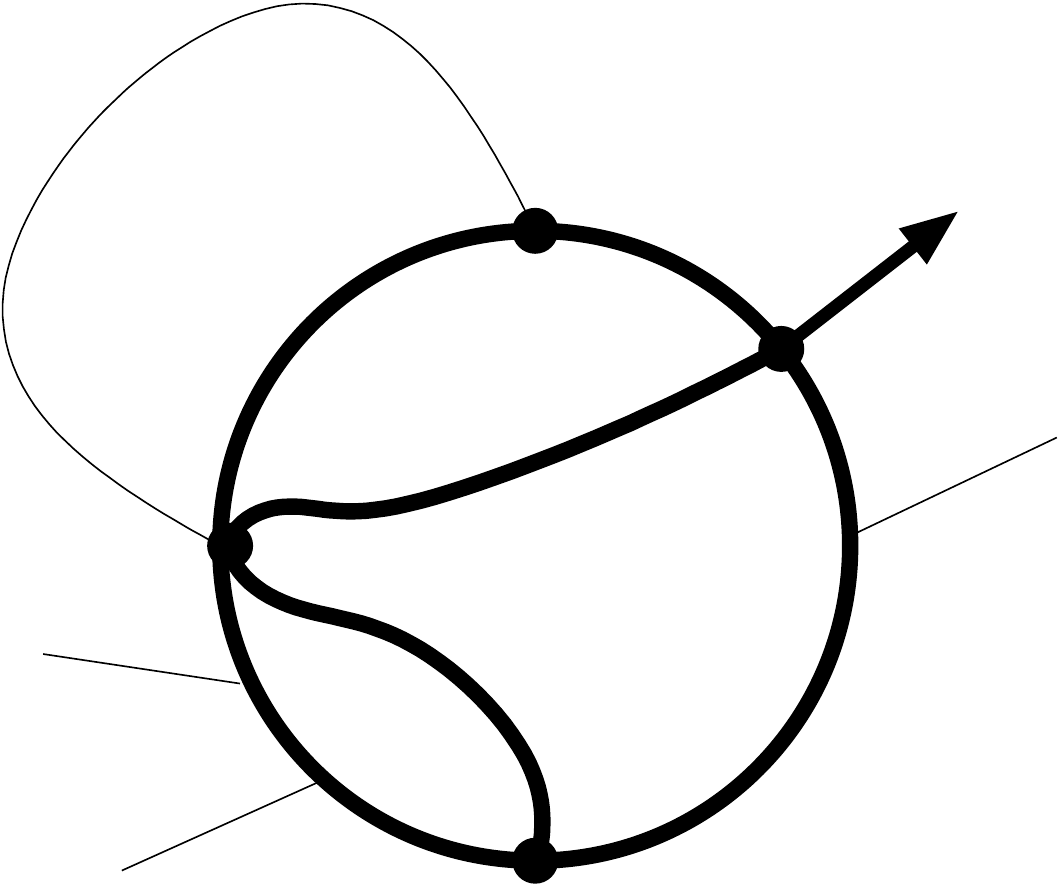}
 $$ 

After the last step has been made, the remaining arrow is linked to the root vertex.
The algorithm is illustrated below for the Kreweras walk $aabbccbac$. We denote the final step by $f$.
$$
\includegraphics[width=11cm]{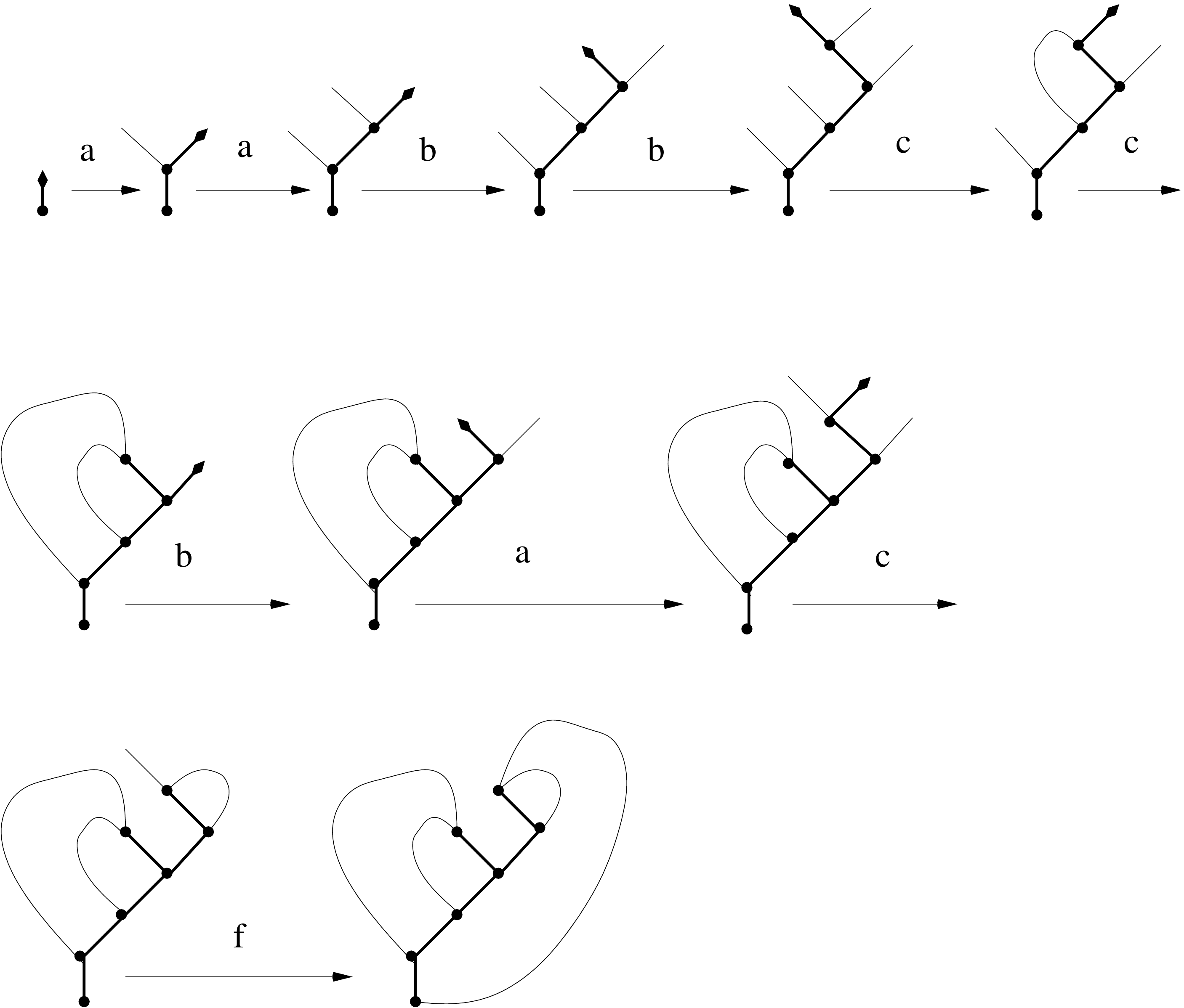}
 $$

It is now a simple matter to check that Bernardi's bijection corresponds exactly to  our 
construction. Indeed by going up in the diagram associated to  the walk we see that application of $\varphi_a$ corresponds to a step of type $a$, similarly for $\varphi_b$. For steps of type $c$ we have to check that our rule is the same as for Bernardi's map $\varphi_c$. This follows from the fact proved by Bernardi that the two half-edges candidates for a pairing with the arrow are ordered in the tree: then the one that is closest to the root must be the one that has lowest numbering in the walk. 

\subsection{Tandem walks, prographs and maps with a bipolar orientation}\label{tandem}
\subsubsection{A bijection between prographs and tandem walks}
\begin{figure}
\caption{From prographs to tandem walks.\newline
(a) A prograph.\newline (b) Cutting left inputs and exploring the tree.\newline (c) The resulting path, recovering the prograph as the dual of the triangulation.}
\includegraphics[height=6cm]{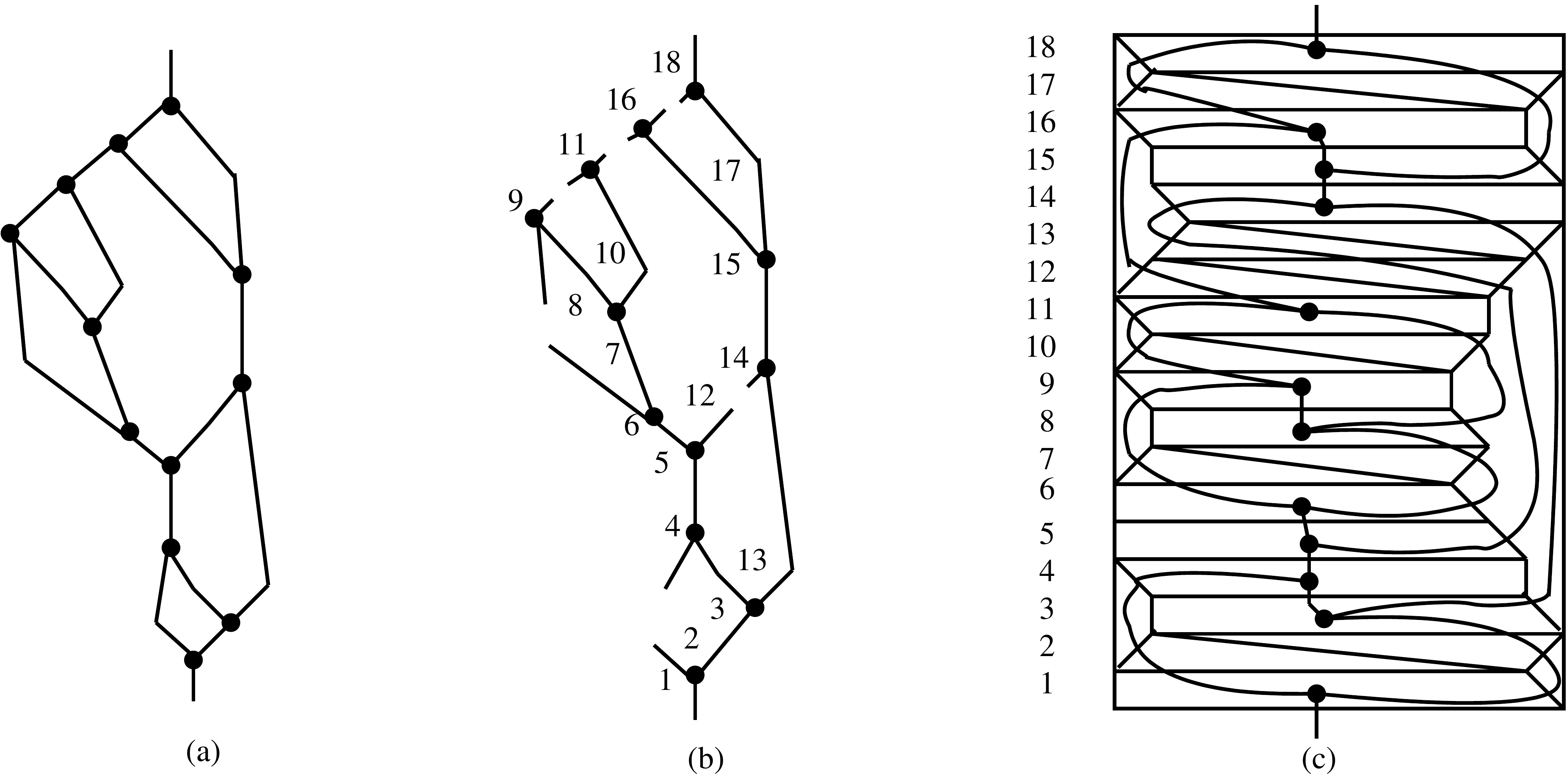}
\label{prograph1}
\end{figure}

We consider oriented planar graphs  
made of   ``products'' with two inputs and one output (from bottom to top)
$${\tt    \setlength{\unitlength}{1pt}
\begin{picture}(100,60)
\thinlines    
\

\put(30,30){\vector(0,1){30}}
\put(10,10){\vector(1,1){18}}
\put(50,10){\vector(-1,1){18}}
\put(27,27){\text{$\bullet$}}

\end{picture}}
$$
and ``coproducts'' with one input and two outputs
$${\tt    \setlength{\unitlength}{1pt}
\begin{picture}(100,60)
\thinlines    
\

\put(30,00){\vector(0,1){30}}
\put(30,30){\vector(1,1){20}}
\put(30,30){\vector(-1,1){20}}
\put(27,27){\text{$\bullet$}}

\end{picture}}
$$
The terminology comes from the theory of operads, cf Borie \cite{Bo}.
One can use outputs as new inputs and get in this way configurations which are oriented from bottom to top. Prographs are such planar configurations with one input and one output. For example here is the only configuration with  one product and one coproduct. 
A more complex example is given in Figure \ref{prograph1}(a).

$${\tt    \setlength{\unitlength}{1pt}
\begin{picture}(100,100)
\thinlines    
\

\put(30,70){\line(0,1){30}}
\put(30,70){\line(1,-1){20}}
\put(30,70){\line(-1,-1){20}}
\put(27,67){\text{$\bullet$}}
\put(30,30){\line(0,-1){30}}
\put(30,30){\line(1,1){20}}
\put(30,30){\line(-1,1){20}}
\put(27,27){\text{$\bullet$}}

\end{picture}}
$$
One can connect the input and the output to $\infty$ and one obtains in this way a planar cubic graph, dual to a triangulation of the sphere. 
Prographs can be enumerated (see \cite{Bo}), they are equinumerous with tandem walks i.e. walks in the quarter plane, starting and ending at $0$,  with steps in the set
$a=(0,1),\ b=(1,-1),\ c=(-1,0)$,
$${\tt    \setlength{\unitlength}{1.5pt}
\begin{picture}(30,35)
\thicklines    

\put(-16,18){$c$}
\put(0,20){\vector(-1,0){10}}

\put(0,20){\vector(1,-1){10}}
\put(13,8){$b$}

\put(0,20){\vector(0,1){10}}
\put(-2,32){$a$}
\end{picture}}
$$ A tandem walk is a word in the alphabet $a,b,c$ and  the condition to stay in the quarter plane is that 
$\sharp a\geq \sharp b\geq \sharp c$ for any prefix. This corresponds to the natural two dimensional generalization of the ballot problem.
 There is a  
simple bijection going from tandem walks in the quarter plane, starting and ending at $0$, to the set of standard Young tableaux with rectangular shape, having three rows of the same length, the numbers in the first line being the positions of the $a$ letters in the word, of the $b$ letters in the second line and of the $c$ letters in the third line. They can thus be counted using the hook formula: there are 
$\frac{2(3n)!}{n!(n+1)!(n+2)!}$ tandem walks, starting and ending at $0$, with $3n$ steps.
Tandem walks can be realized as a subset of $rY$-walks. Indeed, if in a tandem walk we replace each occurence of the step $c=(-1,0)$ by a step $(0,1)$ immediately followed by a step $(-1,-1)$, we obtain an $rY$-walk. Thus we can refer to tandem walks as  $rY$ walks with some forbidden patterns : no step $(0,1)$ is followed by a step $(0,1)$ or a step $(-1,1)$.

Consider a tandem walk in the quarter plane, starting and ending at $0$.
Using the construction of section \ref{sec:3} there are three types of cells:

$${\tt    \setlength{\unitlength}{1pt}
\begin{picture}(300,25)
\thicklines    
\

\put(0,0){\line(1,0){60}}
\put(60,0){\color{blue}\line(-1,1){10}}
\put(0,10){\line(1,0){50}}

\put(0,0){\color{red}\line(0,1){10}}           
\put(20,20){type $a$}

\put(110,10){\line(1,0){50}}
\put(160,10){\color{blue}\line(-1,-1){10}}
\put(100,0){\line(1,0){50}}

\put(100,0){\color{red}\line(1,1){10}}           
\put(110,20){type $b$}

\put(200,10){\line(1,0){60}}
\put(200,10){\color{red}\line(1,-1){10}}
\put(210,0){\line(1,0){50}}

\put(260,0){\color{blue}\line(0,1){10}}           
\put(210,20){type $c$}

\end{picture}}
$$
Quadrilaterals such as  $a$ and $c$  give triangles, while, in case $b$ we will smash the quadrilateral according to its natural bent:
\begin{contraction}
\label{CRtandem}
$${\tt    \setlength{\unitlength}{1pt}
\begin{picture}(300,20)
\thicklines    
\
\put(-10,-6){$w$}\put(57,-6){$z$}\put(-2,13){$u$}\put(59,13){$v$}

\put(0,0){\line(1,0){50}}
\put(60,10){\color{blue}\line(-1,-1){10}}
\put(10,10){\line(1,0){50}}

\put(0,0){\color{red}\line(1,1){10}}

 \put(80,5){$\to$}
\put(100,5){$w$}
\put(110,5){\color{red}\line(1,0){40}}\put(150,5){\color{blue}\line(1,0){40}}\put(148,8){$u$}\put(148,-4){$z$}\put(193,5){$v$}
\end{picture}}
$$
\end{contraction}
\medskip

These rules form  a subset of the Contraction Rules \ref{contract:1}, moreover the other Contraction Rule for $rY$-walks is compatible with the interpretation of a step $c$ as formed by a step $(0,1)$ followed by a step $(-1,-1)$. The left hand of the following picture shows the steps 
$(0,1)$ and $(-1,-1)$ and the right hand shows the $c$ step. Clearly the two are equivalent.

$$
\includegraphics[height=3cm]{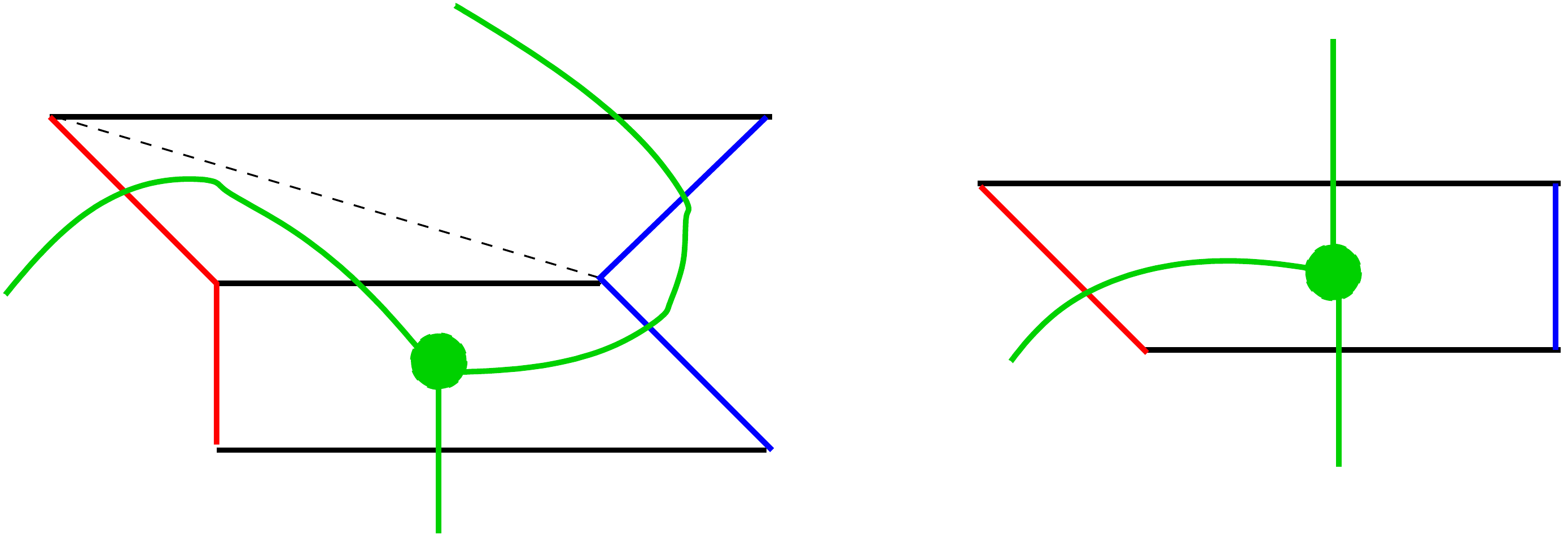}
  $$
The construction using Contraction Rules \ref{CRtandem} associates to a tandem walk in the quarter plane, starting and ending at $0$, a prograph, obtained as the dual of the triangulation.
Indeed it is easy to see that in the construction each $a$ step corresponds to a coproduct while a $c$-step corresponds to a product.
Here is an  example with the word  $abacbc$:
$$
\includegraphics[height=5cm]{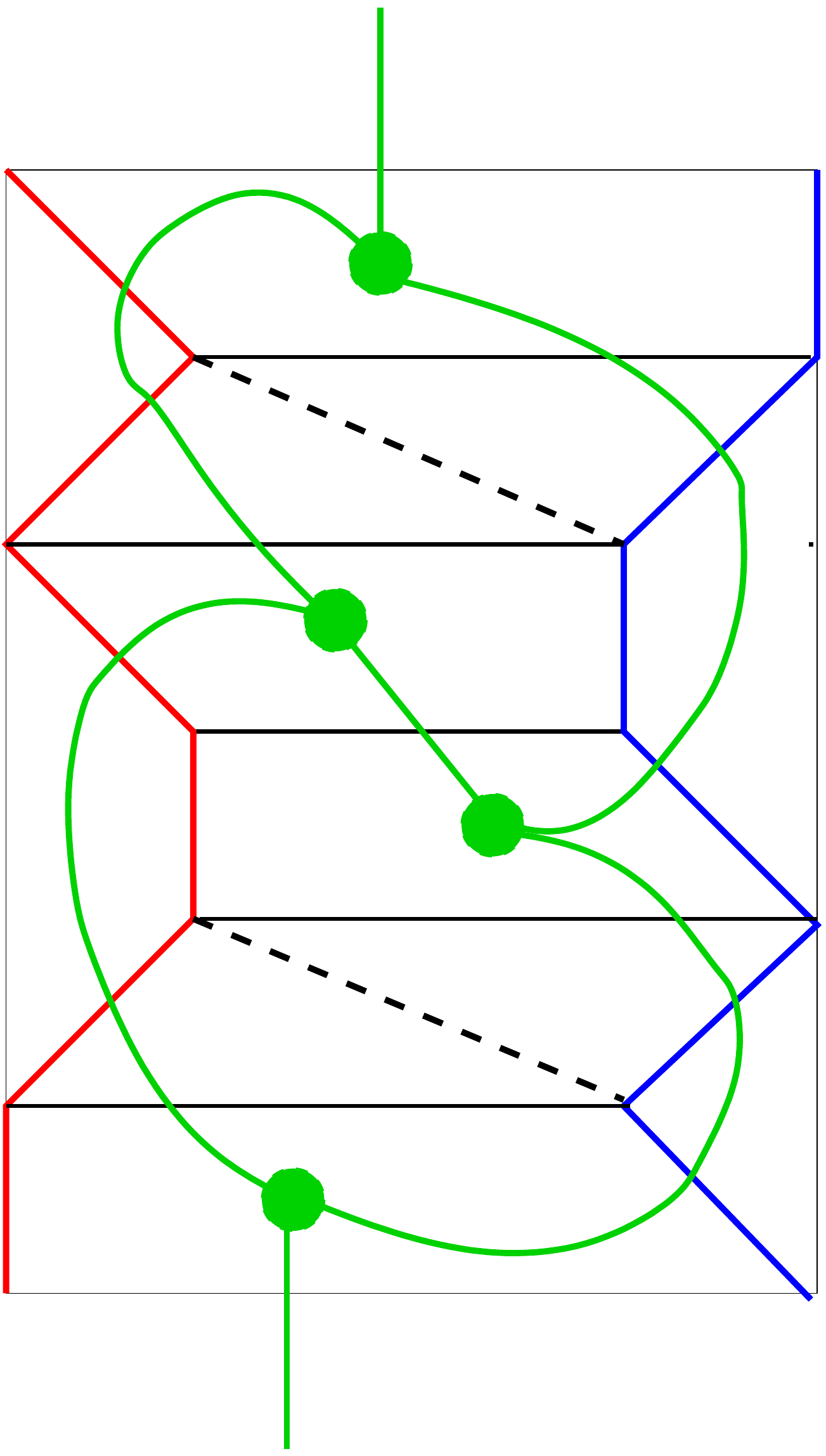}
  $$
The corresponding Young tableau is 
$$\begin{tabular}{|l|r|}
  \hline
  4 & 6 \\
  \hline
  2 & 5 \\
\hline
  1 & 3 \\
  \hline
\end{tabular}$$ 

This correspondance between tandem walks and prographs is in fact a bijection and 
the inverse bijection is readily obtained from the connection, referred to above,  with $rY$-walks. It can be explicitly described as follows. Given a prograph, cut  the left input of each coproduct. The  resulting  graph is a complete spanning tree.   Let us call a corner the region just above a coproduct
(between its two outputs). Make a depth first search of
the tree and order the vertices and corners according to their first appearance. Make an $a$-step the first time you go through  a coproduct, a $b$-step the first time you go through  a product,  and a $c$-step for each corner, this gives the tandem walk corresponding to the prograph.
The construction is shown in Figure \ref{prograph1}.

\subsubsection{Connections with bipolar maps and the bijection of Kenyon, Miller, Sheffield and Wilson }
As in section \ref{sec:3:2}, one can extend the previous construction to walks which do not start or end at zero and to walks with more general step set. Here we will consider walks starting at a point of the type $(0,n)$ and ending  at $(m,0)$ for some $m,n\geq 0$. We will also consider the infinite step set $$\mathcal S=\{(1,-1),(-i,j); i,j\geq 0, i+j>0\}$$ 
The  contruction using Contraction Rules \ref{CRtandem}  provides a map from walks with step set $\mathcal S$ to some set of ``generalized prographs'', which are made of vertices with $i$ inputs and $j$ outputs, as shown below, with 3 inputs and 2 outputs:
$${\tt    \setlength{\unitlength}{1pt}
\begin{picture}(100,60)
\thinlines    
\

\put(40,0){\vector(-1,1){20}}
\put(0,0){\vector(1,1){20}}
\put(20,0){\vector(0,1){20}}

\put(20,20){\vector(1,1){20}}
\put(20,20){\vector(-1,1){20}}
\put(17,17){\text{$\bullet$}}

\end{picture}}
$$
 Instead of looking at these generalized prographs and trying to characterize the ones that we obtain, we will instead look at the dual maps. 
A vertex of a generalized prograph corresponds to a face in the dual map.
The face associated with a step of the form $(-i,j)$ will have $i+j+2$ sides.  We will orient the sides of such a face from west to east, e.g. for $i=3,j=2$ :

$${\tt    \setlength{\unitlength}{1.5pt}
\begin{picture}(100,50)
\thicklines    
\put(30,0){\vector(1,0){50}}
\put(20,10){\vector(1,-1){10}}
\put(10,20){\vector(1,-1){10}}
\put(0,30){\vector(1,-1){10}}
\put(0,30){\line(0,1){20}}
\put(0,50){\vector(1,0){60}}
\put(60,50){\vector(1,-1){10}}
\put(70,40){\vector(1,-1){10}}
\put(80,30){\line(0,-1){30}}
\end{picture}}
$$
For a quadrilateral corresponding to a step $(-1,1)$ :
$$
\includegraphics[height=.5cm]{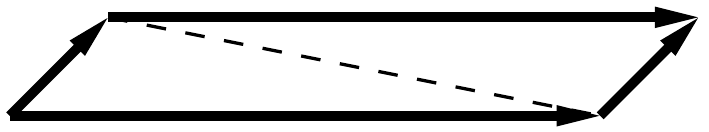}
  $$
Clearly these orientations are compatible with the identifications between sides made when contracting the Motzkin paths or the quadrilaterals. The  orientation of the map produced in this way is acyclic, has a unique source (the west-most point) and a unique sink (the east-most point).
We claim that the oriented maps that we obtain are exactly the  bipolar oriented maps as defined in \cite{KMSW} and,  up to some minor twists, our construction  recovers their   bijection with walks in the quarter plane. 
In \cite{KMSW} this bijection proceeds inductively by looking at the steps of a walk and building an associated map by sewing faces. More precisely a face with $i+j+2$ sides is associated with every step of the form $(-i,j)$. It is easy to check that the sewing algorithm corresponds  exactly to our construction. 
Instead of giving full proofs we will just check the example of \cite{KMSW}, section 2.2 and leave the details to the reader.

\begin{figure}
\caption{(a) The example of \cite{KMSW} with steps\newline $\scriptstyle(1,-1);(0,2);(-1,0);(0,1);(1,-1);(1,-1);(-1,1);(0,1);(1,-1);(1,-1);(1,-1);(1,-1);(-1,0);(-2,1);(1,-1)$.
\newline(b) After a reflection through the main diagonal.
\newline (c) The map drawn as in \cite{KMSW}.
\medskip}
\label{progKMSW}
\includegraphics[width=15cm]{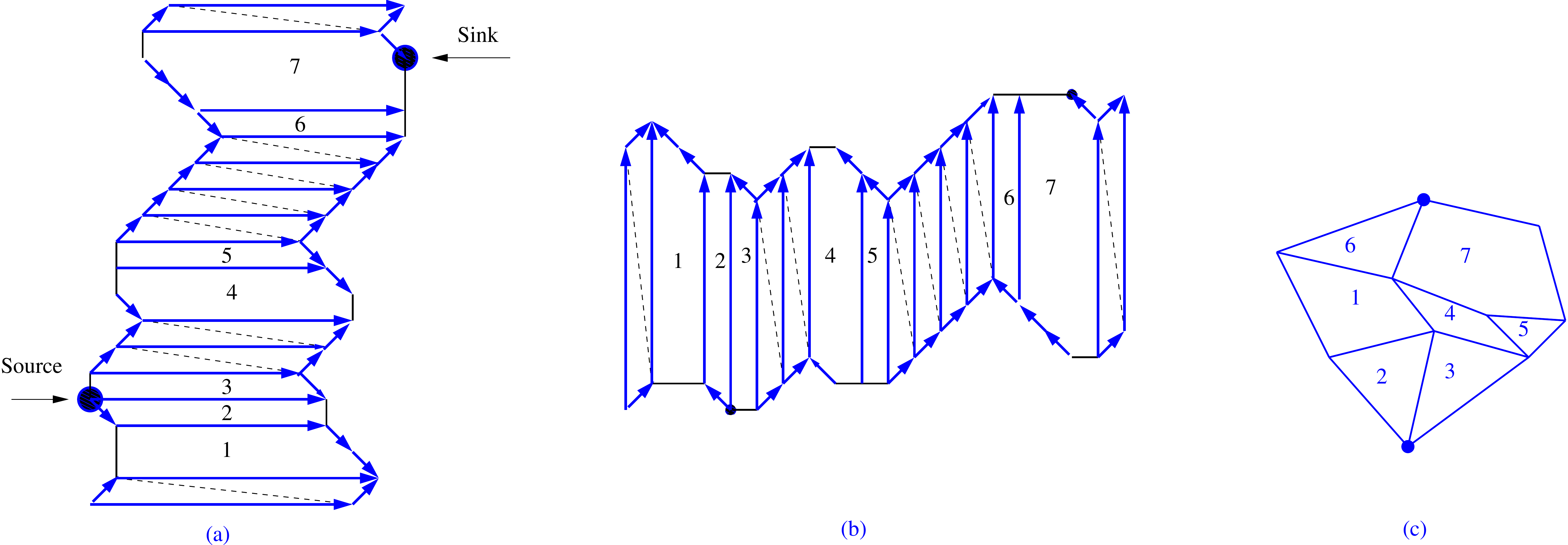}
\end{figure}

In Figure \ref{progKMSW}(a) the picture corresponding to this example  is drawn. 
In order to recover the map of \cite{KMSW}, Figure 4, we make a reflection with respect to the $x=y$ axis, as shown in Figure \ref{progKMSW}(b). After contractions the result  is shown in Figure \ref{progKMSW}(c).

\subsection{Schnyder woods}
We now describe a bijection betwen walks and  Schnyder woods, originated in Li, Sun and Watson \cite{Li-Sun-Watson}.
A Schnyder wood is a planar triangulation in which the three vertices of the external face are coloured,  in clockwise  order, in green, red and blue. The internal edges are also coloured so that they form three trees, one of each colour, rooted on the external vertex of its own colour and containing all internal vertices. These trees are  oriented towards their roots.
The edges have to satisfy the Schnyder condition at each internal vertex: in clockwise order around the vertex, we have successively the  outgoing blue edge, incoming red edges, outgoing green edge, incoming blue edges, outgoing red edge and incoming green edges.

$$
\includegraphics[width=4cm]{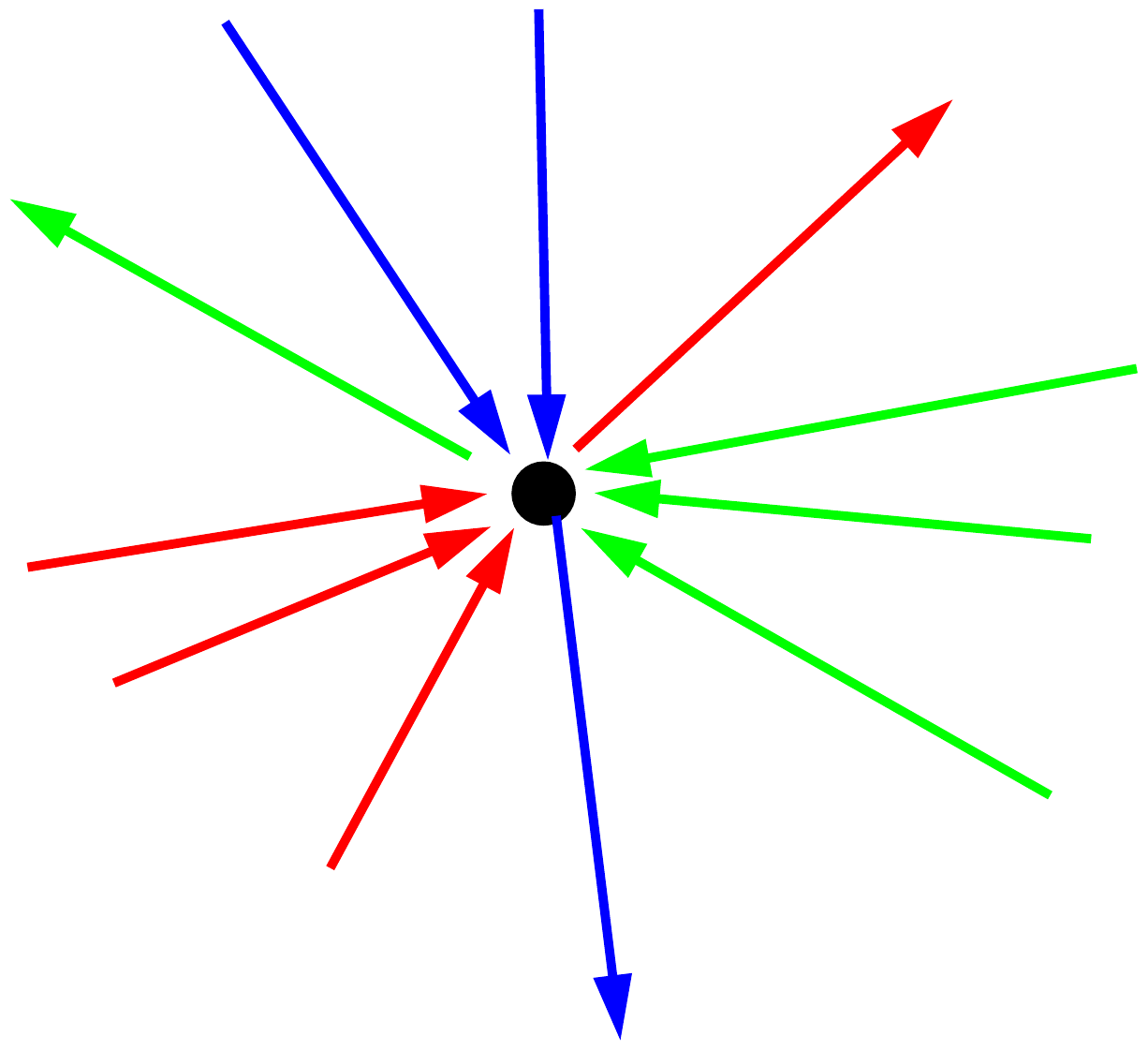}
 $$
It will be convenient for us to also colour in blue the external edge between the blue and red vertices.
 Here I take as running  example the same as Bernardi and Bonichon \cite{BB}, see Figure \ref{Schnyder:1} below.
\begin{figure}

\includegraphics[width=9cm]{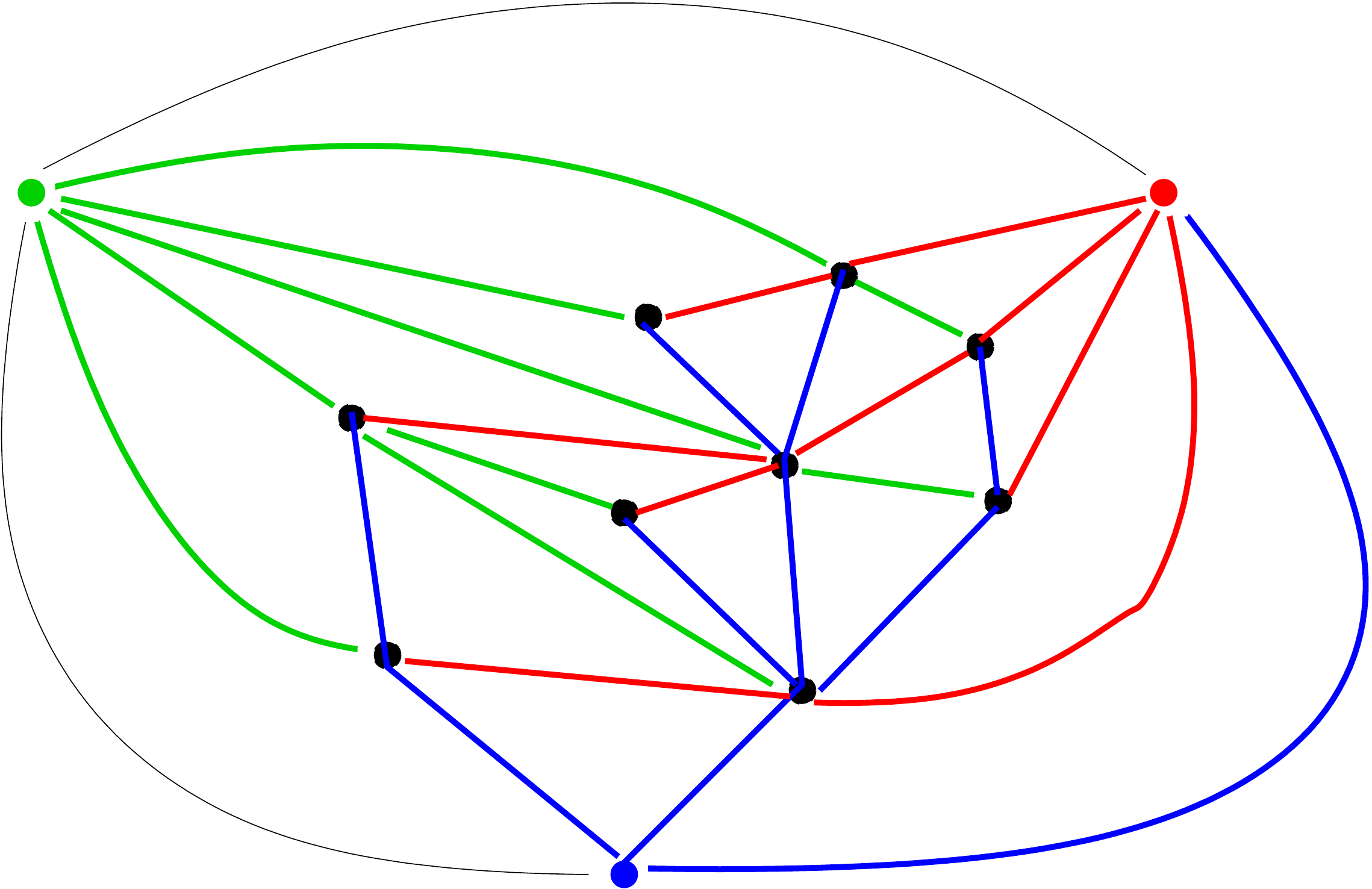}

\caption{A Schnyder wood.}
\label{Schnyder:1}
\end{figure}
 To such a Schnyder wood I will associate a tandem walk in the quarter plane, closely related to the one described in \cite{Li-Sun-Watson}.
 Since the bijection is studied in extensively in \cite{Li-Sun-Watson}, I will only sketch the proof and leave the details to the reader.

We start from a Schnyder wood and construct a tandem walk. Recall that the steps of a tandem walk belong to the set $\{a=(1,0),\ b=(-1,1),\ c=(0,-1)\}$ and the condition to remain in the quarter plane is that $\sharp a\geq \sharp b\geq\sharp c$ for all prefix.
We make a contour exploration of the blue tree, from left to right, starting from the root.
 The first time we encounter a blue edge we make a $a$-step in the walk and when we go down this blue edge we make a $b$ step.
At each vertex we may cross several red edges. Each time we cross an incoming   red edge  we make a $c$-step in the walk. Finally, after we went up the last blue edge, we do not go back to the origin, so we do not make the last $b$-step. Observe that, by the Schnyder condition, a $c$ step can never occur immediately after a $b$ step.
The word obtained in the running example is:
$$aabbacabaccabacbbaacbbbacccc$$

\begin{figure}

\includegraphics[width=6cm]{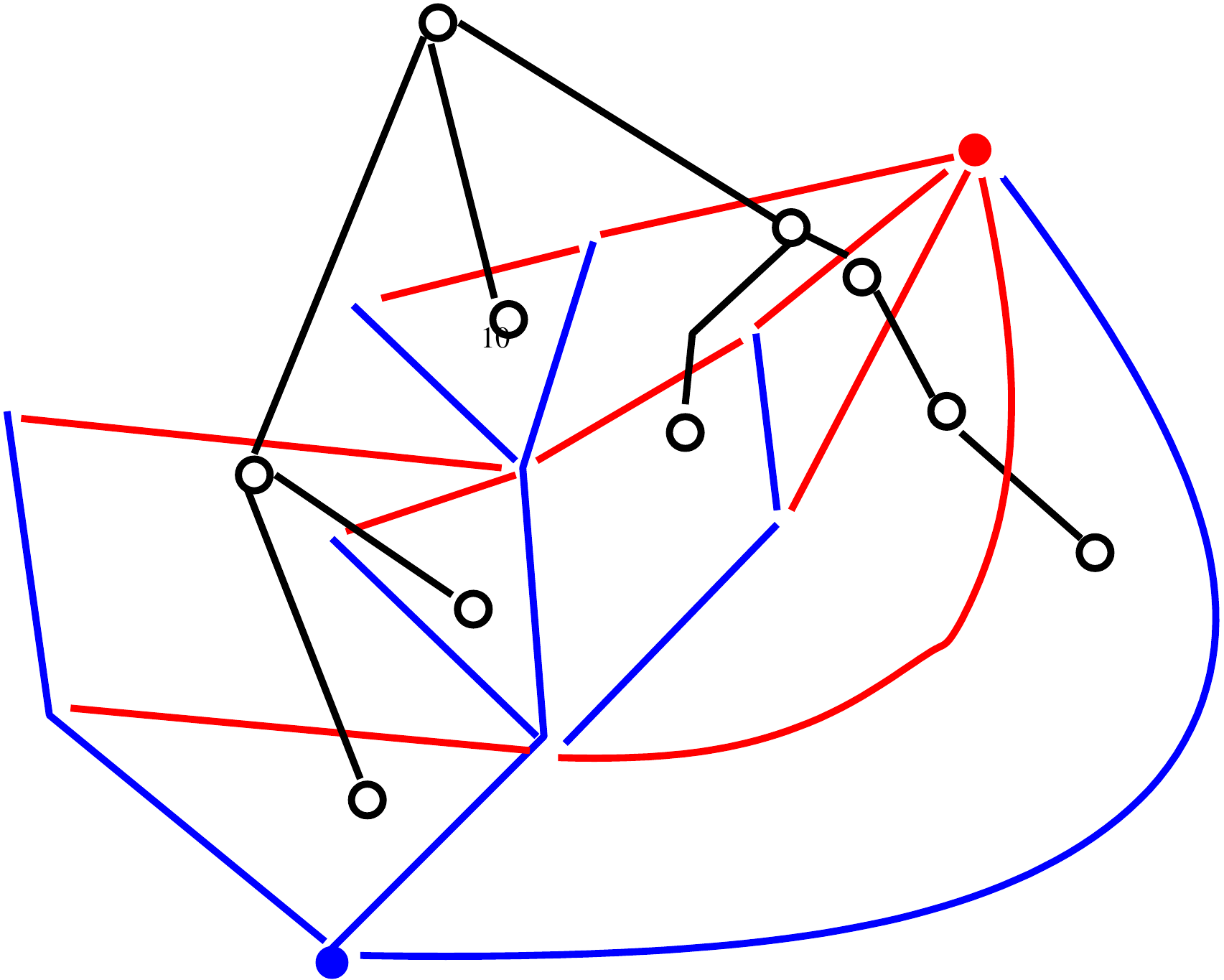}

\caption{A dual tree}
\label{Schnyder:2}
\end{figure}

Let us check that the walk that we obtain is a tandem walk, starting from $(0,0)$ and ending at $(1,0)$. Consider the subword of $a$ and $b$ letters. If we replace the  $a$ and $b$ letters by $u$ and $v$ respectively, we obtain the Dyck path corresponding to the blue tree (with its last step missing), therefore, for any prefix in the word the number of $a$'s is larger than the number of $b$'s. Consider the map obtained by erasing the green edges and the green vertex. In this map construct the tree dual to the red tree and root it on the external face.  It is shown in black in Figure \ref{Schnyder:2}. One can see that the subword formed by the $b$ and $c$ letters gives, upon substituting $u$ for $b$ and $v$ for $c$,  the Dyck path of the exploration of this tree from left to right, therefore, for each prefix in the word, the number of $b$'s is larger than the number of $c$'s. It follows that the word on the letters $a,b,c$ that we obtain from the Schnyder wood gives a tandem walk, from $(0,0)$ to $(1,0)$. 

\begin{figure}
\caption{The Schnyder condition around an internal vertex.}
\label{schnyder-int}
\includegraphics[width=9cm]{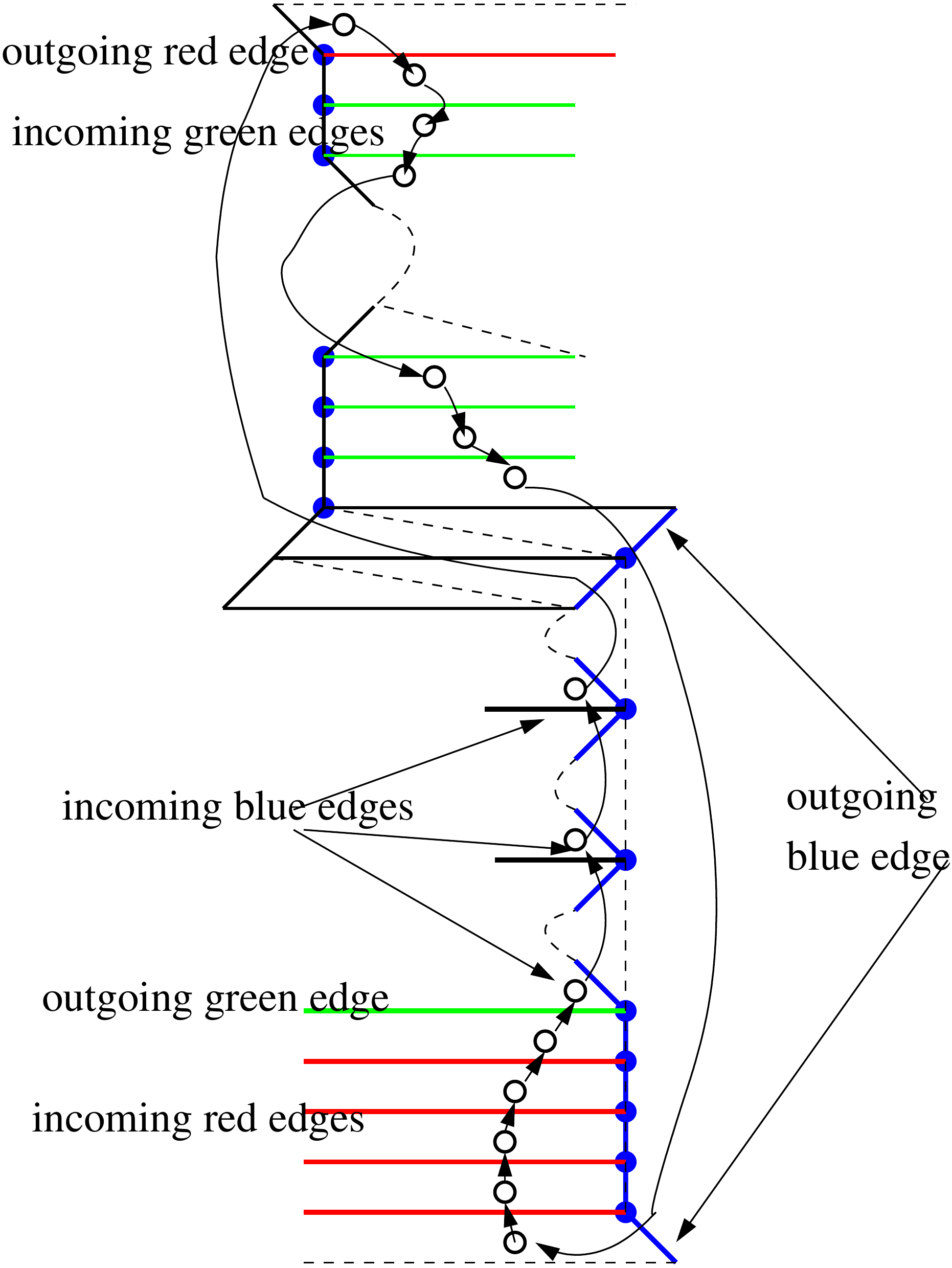}
\end{figure}

Conversely, let $p$ be a tandem walk on the quarter plane from $(0,0)$ to $(1,0)$ satisfying the condition on steps $b$ and $c$ that the pattern $bc$ is forbidden.
We associate to this walk the mating of trees, as in section \ref{sec:2}, using Contraction Rules \ref{CRtandem},
but we do not identify the top and bottom sides of the rectangle. These two sides, together with the last $a$-step, which has not been matched with a $b$-step, form the boundary of the external triangle.
 The blue tree is obtained from the right blue path. It remains to construct the red and green trees. They are obtained by colouring some horizontal edges. We colour red the bottom edge of each triangle of type $c$. Some of the other horizontal edges are sewed to some blue edges and therefore will be coloured in blue. The remaining horizontal edges are coloured in green. We orient the red edges from west to east and the green edges from east to west.
Here is the diagram we obtain in our example. 
$$
\includegraphics[width=5cm]{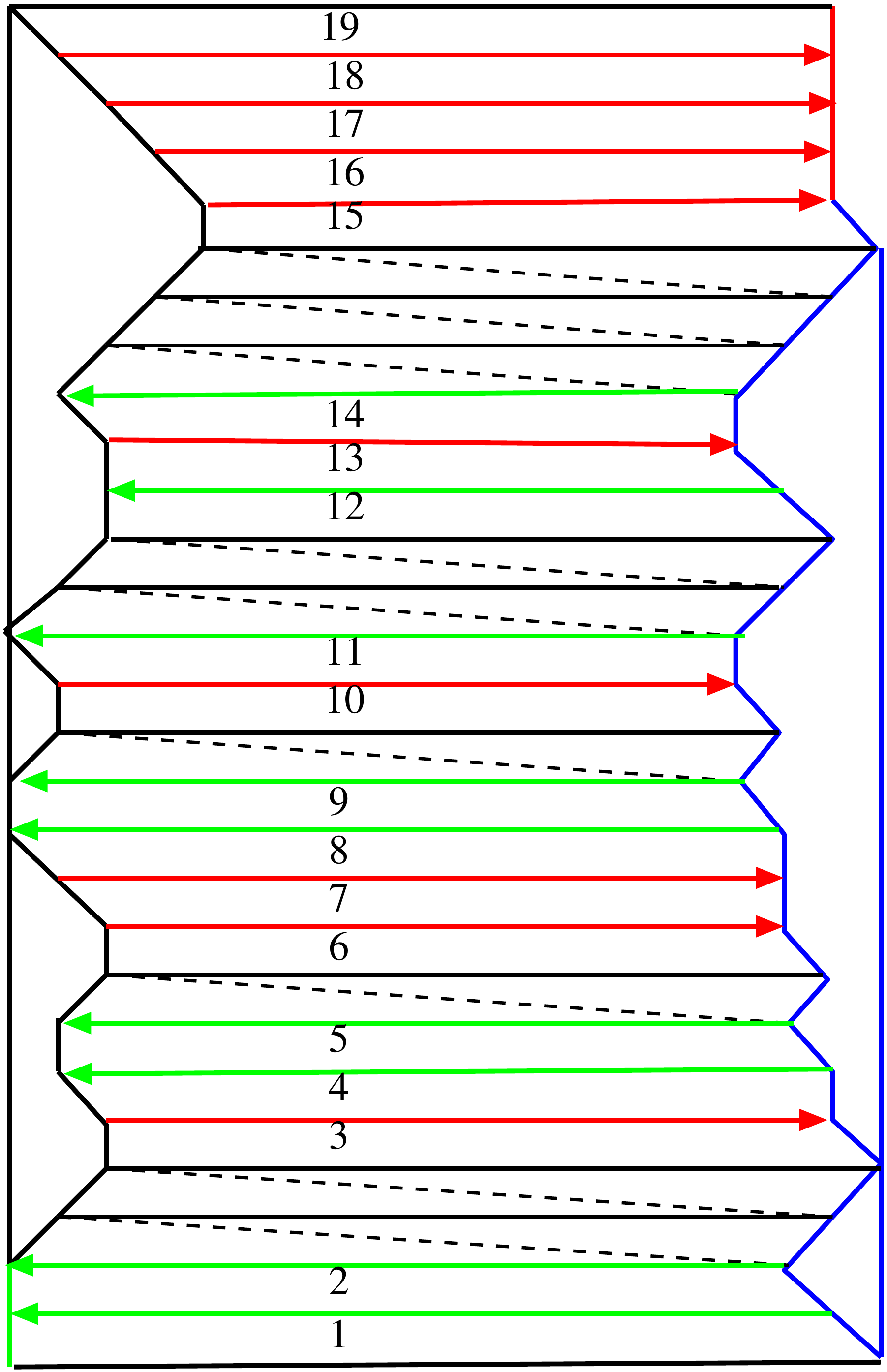}
 $$
We have to check that the coloured triangulation that we have constructed is a Schnyder wood.
 We can identify the vertices of the triangulation, namely they are either the vertices of the blue tree, the green vertex corresponding to the left side of the picture and the red vertex corresponding to the last $c$ steps.
The three external edges are the blue edge from the path, corresponding to the last $a$ step and the upper and lower sides of the picture, as depicted below: 
$${\tt    \setlength{\unitlength}{1.5pt}
\begin{picture}(50,90)
\thicklines    
\put(0,10){\line(1,0){40}}
\put(0,10){\color{green}\line(0,1){70}}
\put(40,10){\color{blue}\line(0,1){30}}
\put(40,40){\color{blue}\line(-1,1){10}}
\put(30,50){\color{red}\line(0,1){30}}
\put(0,80){\line(1,0){30}}
\put(-30,40){green vertex}
\put(20,20){blue vertex}
\put(-5,0){external side}
\put(-5,85){external side}
\put(45,45){blue external side}
\put(25,65){red vertex}

\end{picture}}
$$
The vertices of the left tree which are not the root are all matched to a vertex of the blue tree by the contraction of a quadrangle.
It remains to check the Schnyder condition at each internal vertex. In figure \ref{schnyder-int} we consider a typical  such  vertex.
The blue points  all correspond to the vertex and the successive edges and faces traversed when going clockwise around the vertex are shown on the path with arrows. It is easy to check, using the fact that the pattern $bc$ is forbidden, that these edges follow the Schnyder condition.
Finally, here is the picture of the original Schnyder wood, with the faces numbered as above.

$$
\includegraphics[width=9cm]{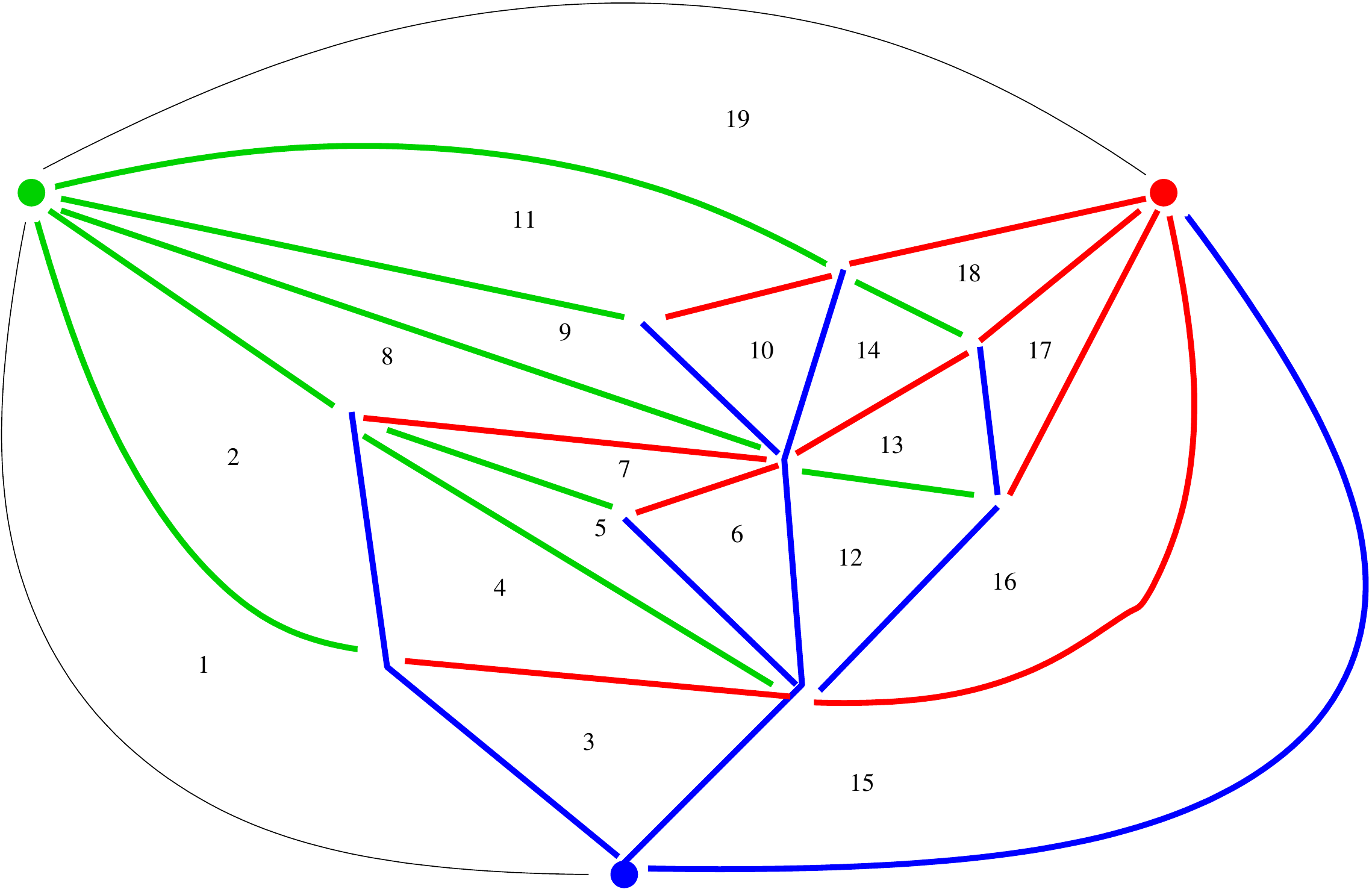}
 $$

\end{document}